\documentclass[10.5pt]{article}
\usepackage{amsfonts}
\usepackage{graphicx}

\usepackage{amsmath}
\usepackage{amssymb}
\usepackage{latexsym}
\usepackage{amsmath, amsfonts,amssymb, amsthm, euscript,makeidx,color,mathrsfs}

\renewcommand{\baselinestretch}{1.2}
\def\sqr#1#2{{\vcenter{\vbox{\hrule height.#2pt
              \hbox{\vrule width.#2pt height#1pt \kern#1pt \vrule width.#2pt}
              \hrule height.#2pt}}}}
\def\signed #1{{\unskip\nobreak\hfil\penalty50
              \hskip2em\hbox{}\nobreak\hfil#1
              \parfillskip=0pt \finalhyphendemerits=0 \par}}
\def\endpf{\signed {$\sqr69$}}

\def\3n{\negthinspace \negthinspace \negthinspace }
\def\2n{\negthinspace \negthinspace }
\def\1n{\negthinspace }

\def\dbE{\mathbb{E}}
\def\dbF{\mathbb{F}}

\def\dbP{\mathbb{P}}

\def\dbR{\mathbb{R}}

\def\sX{\mathscr{X}}
\def\sY{\mathscr{Y}}

\def\={\buildrel \triangle \over =}

\def\ds{\displaystyle}

\def\ns{\noalign{\ss}}
\def\a{\alpha}

\def\d{\delta}

\def\z{\zeta}

\def\l{\lambda}
\def\m{\mu}

\def\si{\sigma}
\def\t{\tau}
\def\f{\varphi}
\def\th{\theta}

\def\G{\Gamma}
\def\D{\Delta}
\def\Th{\Theta}
\def\L{\Lambda}

\def\O{\Omega}

\def\cF{{\cal F}}

\def\cL{{\cal L}}

\def\cP{{\cal P}}

\def\cS{{\cal S}}

\def\cY{{\cal Y}}

\def\ss{\smallskip}
\def\ms{\medskip}
\def\bs{\bigskip}
\def\q{\quad}
\def\qq{\qquad}
\def\hb{\hbox}

\def\lan{\mathop{\langle}}
\def\ran{\mathop{\rangle}}

\def\h{\widehat}
\def\wt{\widetilde}

\def\cd{\cdot}
\def\cds{\cdots}

\def\as{\hbox{\rm a.s.}}

\def\les{\leqslant}
\def\ges{\geqslant}

\def\({\Big (}
\def\){\Big )}
\def\[{\Big[}
\def\]{\Big]}
\def\bde{\begin{definition}\label}
\def\ede{\end{definition}}
\def\be{\begin{equation}}
\def\bel{\begin{equation}\label}
\def\ee{\end{equation}}
\def\bt{\begin{theorem}\label}
\def\et{\end{theorem}}
\def\bc{\begin{corollary}\label}
\def\ec{\end{corollary}}
\def\bl{\begin{lemma}\label}
\def\el{\end{lemma}}
\def\bp{\begin{proposition}\label}
\def\ep{\end{proposition}}
\def\bas{\begin{assumption}}
\def\eas{\end{assumption}}
\def\br{\begin{remark}\label}
\def\er{\end{remark}}
\def\ba{\begin{array}}
\def\ea{\end{array}}

\def\rf{\eqref}

\def\square#1{\vbox{\hrule\hbox{\vrule height#1%
     \kern#1\vrule}\hrule}}
\def\rectangle#1#2{\vbox{\hrule\hbox{\vrule height#1%
     \kern#2\vrule}\hrule}}

\font\tenbb=msbm10 \font\sevenbb=msbm7 \font\fivebb=msbm5

\newfam\bbfam
\scriptscriptfont\bbfam=\fivebb \textfont\bbfam=\tenbb
\scriptfont\bbfam=\sevenbb

\newtheorem{theorem}{Theorem}[section]
\newtheorem{corollary}[theorem]{Corollary}

\newtheorem{lemma}[theorem]{Lemma}
\newtheorem{proposition}[theorem]{Proposition}

\theoremstyle{definition}
\newtheorem{definition}[theorem]{Definition}
\newtheorem{remark}[theorem]{Remark}

\makeatletter
   
   \@addtoreset{equation}{section}
\makeatother

\begin{document}

\title{\bf Backward Stochastic Volterra Integral Equations \\
--- Representation of Adapted Solutions \footnote{This work is supported in part by NSF of China (Grant 11401404, 11471231) and NSF Grant
DMS-1406776.}}

\author{Tianxiao Wang\footnote{School of Mathematics, Sichuan University, Chengdu 610065, China}~~~and~~Jiongmin Yong\footnote{Department of
Mathematics, University of Central Florida, Orlando, FL 32816, USA}}

\maketitle

\begin{abstract} For backward stochastic Volterra integral equations (BSVIEs, for short), under some mild conditions, the so-called adapted solutions or adapted M-solutions uniquely exist. However, satisfactory regularity of the solutions is difficult to obtain in general. Inspired by the decoupling idea of forward-backward stochastic differential equations, in this paper, for a class of BSVIEs, a representation of adapted M-solutions is established by means of the so-called representation partial differential equations and (forward) stochastic differential equations. Well-posedness of the representation partial differential equations are also proved in certain sense.

\end{abstract}

\ms


\section{Introduction.}\label{1}

Let $(\O,\cF,\dbF,\dbP)$ be a complete filtered probability space on which a standard $d$-dimensional Brownian motion $W(\cd)$ is defined, with $\dbF\equiv\{\cF_t\}_{t\ges0}$ being its natural filtration augmented by all the $\dbP$-null sets. We consider the following stochastic integral equation in $\dbR^m$:
\bel{BSVIE2}Y(t)=\psi(t)+\int_t^Tg(t,s,Y(s),Z(t,s),Z(s,t))ds-\int_t^TZ(t,s)dW(s),\qq t\in[0,T].\ee
Such an equation is called a {\it backward stochastic Volterra integral equation} (BSVIE, for short). In the above, $\psi(\cd)$, called a {\it free term}, is an $\cF_T$-measurable stochastic process (not necessarily $\dbF$-adapted) and $g(\cd)$, called the {\it generator} of the above BSVIE, is a given map, deterministic or random. The unknown that we are looking for is the pair $(Y(\cd),Z(\cd\,,\cd))$. Let us look at a special case of the above BSVIE. Suppose
$$g(t,s,y,z,\z)=g(s,y,z),\qq\psi(t)=\xi,\qq\forall(t,s,y,z,\z),$$
with $\xi$ being an $\cF_T$-measurable random variable, and $g(\cd)$ being a proper map.
Then the above BSVIE is reduced to the following form:
\bel{BSVIE1*}Y(t)=\xi+\int_t^Tg(s,Y(s),Z(t,s))ds-\int_t^TZ(t,s)dW(s),\qq t\in[0,T].\ee
It is comparable with the integral form of {\it backward stochastic differential equation} (BSDE, for short) which takes the following form:
\bel{BSDE1}Y(t)=\xi+\int_t^Tg(s,Y(s),Z(s))ds-\int_t^TZ(s)dW(s),\qq t\in[0,T].\ee
Now, if BSDE \eqref{BSDE1} admits a unique {\it adapted solution} $(Y(\cd),Z(\cd))$, by which we mean that this pair is $\dbF$-adapted and satisfies (\ref{BSDE1}) in the usual It\^o sense, then, this solution must also be an {\it adapted solution} of BSVIE \eqref{BSVIE1*} with $Z(t,s)\equiv Z(s)$. From this point of view, BSVIE can be regarded as an extension of BSDE.

\ms

Linear BSDEs were firstly introduced by Bismut in 1973 (\cite{Bismut 1973}) while he was studying stochastic linear-quadratic optimal control problems. In 1990, Pardoux--Peng generalized Bismut's linear BSDEs to general nonlinear BSDEs (\cite{Pardoux-Peng 1990}). Shortly after, BSDE theory was found to have very interesting applications in mathematical finance (see for example, \cite{Duffie-Epstein 1992, El Karoui-Peng-Quenez 1997}), and many other developments have been appearing thereafter, including nonlinear Feynman-Kac formula, nonlinear expectations, dynamic risk measures, path-dependent partial differential equations, etc., see for examples \cite{Peng 2004, Peng 2010, Soner-Touzi-Zhang 2012, Ekren-Keller-Touzi-Zhang 2014}, and references cited therein. Relevant theory of BSDEs can also be found in \cite{Briand-Delyon-Hu-Pardoux-Stoica 2003, Kobylanski 2000, Briand-Hu 2006, Delbaen-Hu-Bao 2011}. On the other hand, an extension of BSDE to the so-called {\it forward-backward stochastic differential equations} (FBSDEs, for short) was initiated by Antonelli in 1993 (\cite{Antonelli 1993}). A general form of FBSDE takes the following form:
\bel{FBSDE1}\left\{\2n\ba{ll}
\ds dX(t)=b(t,X(t),Y(t),Z(t))dt+\si(t,X(t),Y(t),Z(t))dW(t),\qq t\in[0,T],\\
\ns\ds dY(t)=-g(t,X(t),Y(t),Z(t))dt+Z(t)dW(t),\qq t\in[0,T],\\
\ns\ds X(0)=x,\qq Y(T)=h(X(T)),\ea\right.\ee
for some maps $b,\si,g,h$. General theories on FBSDEs were developed in the past two and half decays, see \cite{Ma-Protter-Yong 1994, Ma-Yong 1995, Hu-Peng 1995, Yong 1997, Ma-Yong 1999}, and references cited therein. A triple $(X(\cd),Y(\cd),Z(\cd))$ of processes is called an {\it adapted solution} to (\ref{FBSDE1}) if it is $\dbF$-adapted, and satisfies (\ref{FBSDE1}) in the usual It\^o's sense. It is known that under some very general conditions, FBSDE (\ref{FBSDE1}) admits a unique adapted solution $(X(\cd),Y(\cd),Z(\cd))$, and the following estimate holds (for some $p>1$): (\cite{Ma-Yong 1999,Yong 1997,Yong 2010})
\bel{estimate1}\dbE\Big\{\sup_{t\in[0,T]}|X(t)|^p+\sup_{t\in[0,T]}|Y(t)|^p
+\(\int_0^T|Z(t)|^2dt\)^{p\over2}\Big\}\les K\big(1+|x|^p\big),\q\forall x\in\dbR^n.\ee
Hereafter, $K>0$ will stand for a generic constant which can be different from line to line. We should point out that in general, $Z(\cd)$ only belongs to the following space:
$$\ba{ll}
\ns\ds L^p_\dbF(\O;L^2(0,T;\dbR^{m\times d}))=\Big\{Z:[0,T]\times\O\to\dbR^{m\times d}\bigm|t\mapsto Z(t)\hb{ is $\dbF$-adapted, }\dbE\(\int_0^T|Z(t)|^2dt\)^{p\over2}<\infty\Big\}.\ea$$
Therefore, $t\mapsto Z(t)$ is not necessarily continuous and for a given $t\in[0,T]$, $Z(t)$ might not be well-defined. This leads to some difficulties in numerical aspects for adapted solutions to BSDEs, and also restricts the usage of BSDEs/FBSDEs in applications. On the other hand, the so-called {\it Four Step Scheme} (\cite{Ma-Protter-Yong 1994, Ma-Yong 1999}, see also \cite{Hu-Ma 2004,Ma-Yong-Zhao 2010,Ma-Wu-Zhang-Zhang 2015}) for FBSDEs gives a representation of the adapted solution $(Y(\cd),Z(\cd))$ of the BSDE in (\ref{FBSDE1}) via a solution to a relevant partial differential equation (PDE, for short), together with the solution $X(\cd)$ to the (forward) stochastic differential equation (FSDE, for short) in (\ref{FBSDE1}). Such kind of representation can help people to overcome the difficulties encountered in designing numerical algorithms for BSDEs/FBSDEs \cite{Zhang 2004}. This also substantially broadens the applicability of BSDEs/FBSDEs in solving real problems.

\ms

Since the representation of adapted solutions to BSDEs/FBSDEs is very closely related to the main results in the current paper, and the techniques/ideas will be used below, we now elaborate the procedure here (See \cite{Ma-Protter-Yong 1994,Ma-Yong 1999} for more details).

\ms

Consider FBSDE \eqref{FBSDE1}. Suppose all the involved functions $b,\si,g,h$ are deterministic. Inspired by the so-call {\it invariant embedding} (\cite{Bellman-Kalaba-Wing 1960,Bellman-Wing 1975}), one could expect that there is a relation between the backward component $Y(\cd)$ and the forward component $X(\cd)$ as follows:
$$Y(t)=\Th(t,X(t)),\qq t\in[0,T],$$
for some differentiable function $\Th:[0,T]\times\dbR^n\to\dbR^m$. If this is the case, then by It\^o's formula, one must have
$$\ba{ll}
\ns\ds-g(t,X(t),Y(t),Z(t))dt+Z(t)dW(t)=dY(t)=\[\Th_t(t,X(t))+\Th_x(t,X(t))b(t,X(t),Y(t),Z(t))\\
\ns\ds+{1\over2}\,\si(t,X(t),Y(t),Z(t))^\top\Th_{xx}(t,X(t))\si(t,X(t),Y(t),Z(t))\]dt+\Th_x(t,X(t)\si(t,X(t),Y(t),Z(t))dW(t),\ea$$
where
\bel{si Th si}\ba{ll}
\ns\ds\si(t,x,y,z)^\top\Th_{xx}(t,x)\si(t,x,y,z)=\sum_{k=1}^d
\begin{pmatrix}\si_k(t,x,y,z)^\top\Th^1_{xx}(t,x)\si_k(t,x,y,z)\\
\si_k(t,x,y,z)^\top\Th^2_{xx}(t,x)\si_k(t,x,y,z)\\
\vdots\\
\si_k(t,x,y,z)^\top\Th^m_{xx}(t,x)\si_k(t,x,y,z)\end{pmatrix},\ea\ee
with
$$\si(t,x,y,z)=\(\si_1(t,x,y,z),\si_2(t,x,y,z),\cds,\si_d(t,x,y,z)\),\qq\Th(t,x)
=\begin{pmatrix}\Th^1(t,x)\\ \Th^2(t,x)\\ \vdots\\ \Th^m(t,x)\end{pmatrix}.$$
Therefore, the following should hold:
$$Z(t)=\Th_x(t,X(t))\si(t,X(t),Y(t),Z(t)),\qq t\in[0,T].$$
Suppose there exists a well-defined map $\z:[0,T]\times\dbR^n\times\dbR^m\to\dbR^{m\times d}$ satisfying
$$\z(t,x,y)=\Th_x(t,x)\si\big(t,x,y,\z(t,x,y)\big),\qq(t,x,y)\in[0,T]\times\dbR^n\times\dbR^m,$$
which is the case, trivially, if $\si(t,x,y,z)\equiv\si(t,x,y)$ (independent of $z$). Then $\Th(\cd\,,\cd)$ should solve the following system of quasi-linear parabolic PDE system:
\bel{PDE1}\left\{\2n\ba{ll}
\ds\Th_t(t,x)+{1\over2}\,\si\big(t,x,\Th(t,x),\z(t,x,\Th(t,x))\big)^\top\Th_{xx}(t,x)
\si\big(t,x,\Th(t,x),\z(t,x,\Th(t,x))\big)\\
\ns\ds\q+\Th_x(t,x)b\big(t,x,\Th(t,x),\z(t,x,\Th(t,x))\big)
+g\big(t,x,\Th(t,x),\z(t,x,\Th(t,x))\big)=0,\\
\ns\ds\qq\qq\qq\qq\qq\qq\qq\qq\qq\qq\qq\qq(t,x)
\in[0,T]\times\dbR^n,\\
\ns\ds\Th(T,x)=h(x),\qq x\in\dbR^n.\ea\right.\ee
Suppose the above PDE has a classical solution $\Th(\cd\,,\cd)$. Then we solve the following FSDE:
\bel{FSDE}\left\{\2n\ba{ll}
\ds dX(t)=b\big(t,X(t),\Th(t,X(t)),\z(t,X(t),\Th(t,X(t)))\big)dt\\
\ns\ds\qq\qq\q+\si\big(t,X(t),\Th(t,X(t)),\z(t,X(t),\Th(t,X(t)))\big)dW(t),\qq t\in[0,T],\\
\ns\ds X(0)=x.\ea\right.\ee
Now, if the above FSDE admits a solution $X(\cd)$, then we have the representation of the backward components $(Y(\cd),Z(\cd))$ in terms of the forward component $X(\cd)$:
\bel{Rep1}Y(t)=\Th(t,X(t)),\q Z(t)=\z(t,X(t),\Th(t,X(t))),\qq t\in[0,T].\ee
In the above, $\Th(\cd\,,\cd)$ is called a {\it decoupling field} of the FBSDE (\ref{FBSDE1}) (\cite{Ma-Wu-Zhang-Zhang 2015}), and \rf{PDE1} is called the {\it representation PDE} since the solution $\Th(\cd\,,\cd)$ allows us to represent the backward component $(Y(\cd),Z(\cd))$ in terms of the forward component $X(\cd)$. From the above, we see that as long as all the involved functions are nice enough (in a suitable sense), the above representation (\ref{Rep1}) provides useful regularity information on $(Y(\cd),Z(\cd))$, especially for $Z(\cd)$. This actually has played some interesting roles in numerical aspects of BSDEs/FBSDEs (\cite{Douglas-Ma-Protter 1996,Zhang 2004}).

\ms

Note that in the case that both $b$ and $\si$ are independent of $(Y(\cd),Z(\cd))$, for which the FBSDE is decoupled, the representation PDE becomes
\bel{PDE2}\left\{\2n\ba{ll}
\ds\Th_t(t,x)+{1\over2}\,\si(t,x)^\top\Th_{xx}(t,x)
\si(t,x)+\Th_x(t,x)b(t,x)
+g\big(t,x,\Th(t,x),\Th_x(t,x)\si(t,x)\big)=0,\\
\ns\ds\qq\qq\qq\qq\qq\qq\qq\qq\qq\qq\qq\qq\qq\qq(t,x)
\in[0,T]\times\dbR^n,\\
\ns\ds\Th(T,x)=h(x),\qq x\in\dbR^n,\ea\right.\ee
whose solvability conditions are much simpler than those for (\ref{PDE1}). In this case, (\ref{Rep1}) becomes
\bel{Rep2}Y(t)=\Th(t,X(t)),\q Z(t)=\Th_x(t,X(t))\si(t,X(t)),\qq t\in[0,T],\ee
with $X(\cd)$ being the solution of FSDE:
\bel{SDE2}\left\{\2n\ba{ll}
\ds dX(t)=b(t,X(t))dt+\si(t,X(t))dW(t),\qq t\in[0,T],\\
\ns\ds X(0)=x.\ea\right.\ee

Let us return to FBSDE \eqref{FBSDE1}. For any $(s,x)\in[0,T)\times\dbR^n$, let $(X(\cd\,;s,x),Y(\cd\,;s,x),Z(\cd\,;s,x))$ be the (unique) adapted solution to \eqref{FBSDE1} on $[s,T]$ with $X(0)=x$ replaced by $X(s)=x$. Then
\bel{Feynman-Kac}\Th(s,x)=Y(s;s,x),\qq(s,x)\in[0,T)\times\dbR^n.\ee
Thus, the solution $\Th(\cd\,,\cd)$ to the PDE \eqref{PDE1} admits a representation $Y(\cd\,;\cd\,,\cd)$, a part of the adapted solution to FBSDE \eqref{FBSDE1}. This is called a {\it nonlinear Feynman-Kac formula} (see \cite{Peng 1991}, \cite{Ma-Zhang 2002}).

\ms

We now consider BSVIEs. In 2002, Lin firstly introduced a BSVIE (\cite{Lin 2002}) as an extension of BSDEs, in which the term $Z(s,t)$ did not appear and $\psi(t)\equiv\xi$ is a fixed $\cF_T$-measurable random variable. The form (\ref{BSVIE2}), including the term $Z(s,t)$ with general $\psi(\cd)$, was firstly introduced by the second author of the current paper in 2006 (\cite{Yong 2006}), motivated by optimal control of (forward) stochastic Volterra integral equations (FSVIEs, for short). When $Z(s,t)$ is absent, the BSVIE (\ref{BSVIE2}) becomes:
\bel{BSVIE1}Y(t)=\psi(t)+\int_t^T\1n g(t,s,Y(s),Z(t,s))\,ds-\int_t^TZ(t,s)dW(s),\qq t\in[0,T].\ee
Hereafter, we call (\ref{BSVIE1}) and (\ref{BSVIE2}) Type-I and Type-II BSVIEs, respectively. Thus, Type-I BSVIE is a special case of Type-II BSVIE.

\ms

Mimicking the case of BSDEs, a pair $(Y(\cd),Z(\cd\,,\cd))$ is called an {\it adapted solution} to BSVIE (\ref{BSVIE2}) if for each $t\in[0,T)$, the map $s\mapsto(Y(s),Z(t,s))$ is $\dbF$-adapted on $[t,T]$, and satisfies equation (\ref{BSVIE2}) in the usual It\^o sense. For Type-I BSVIE (\ref{BSVIE1}), one needs only to determine $Z(t,s)$ for $(t,s)\in\D[0,T]$, where
\bel{D}\D[0,T]=\Big\{(t,s)\in[0,T]^2\bigm|0\les t\les s\les T\Big\}.\ee
Therefore, under proper conditions, a Type-I BSVIE admits a unique adapted solution. However, for a Type-II BSVIE, due to the appearance of $Z(s,t)$ in the equation, we need to determine $Z(t,s)$ for $(t,s)\in[0,T]^2$, and (\ref{BSVIE2}) alone does not give enough restrictions on $Z(t,s)$. Consequently, as pointed out in \cite{Yong 2008}, the adapted solution to Type-II BSVIE (\ref{BSVIE2}) is not unique. Inspired by the duality principle needed in the optimal control of FSVIEs, the so-called {\it adapted M-solution} was introduced in \cite{Yong 2008}: A pair $(Y(\cd),Z(\cd\,,\cd))$ is called an adapted M-solution to (\ref{BSVIE2}) if it is an adapted solution and moreover, the following holds:
\bel{M}Y(t)=\dbE[Y(t)]+\int_0^tZ(t,s)dW(s),\q t\in[0,T],~\as\ee
It was proved in \cite{Yong 2008} that under certain conditions, Type-II BSVIE (\ref{BSVIE2}) admits a unique adapted M-solution. Moreover, the following estimate holds:
\bel{estimate2}\ba{ll}
\ns\ds\dbE\(\int_0^T|Y(t)|^2ds+\int_0^T\3n\int_0^T|Z(t,s)|^2dsdt\)\les K\(\int_0^T\int_t^T|g(t,s,0,0,0,0)|^2dsdt+\dbE\int_0^T|\psi(s)|^2ds\).\ea\ee
For some relevant results on BSVIEs, the readers are further referred to \cite{Shi-Wang-Yong 2013,Shi-Wany-Yong 2015,Wang-Yong 2015}.

\ms

From \cite{Yong 2008}, we see that to get some further regularities beyond the above estimate \eqref{estimate2} for the process $(Y(\cd),Z(\cd\,,\cd))$, many technical conditions have to be imposed, the proofs are quite technical, and unfortunately, the regularity results were still not satisfactory, especially that of the process $(t,s)\mapsto Z(t,s)$.

\ms

Inspired by the decoupling FBSDEs presented above, we naturally ask: Is it possible to get representation of adapted solutions for Type-I BSVIEs and adapted M-solutions for Type-II BSVIEs similar to (\ref{Rep2}) for BSDEs? More precisely, we will consider the following BSVIEs:
\bel{BSVIE1**}\ba{ll}
\ns\ds Y(t)=\psi(t,X(t),X(T))+\int_t^Tg(t,s,X(t),X(s),Y(s),Z(t,s))ds-\int_t^T
Z(t,s)dW(s),\\
\ns\ds\qq\qq\qq\qq\qq\qq\qq\qq\qq\qq\qq\qq\qq\qq t\in[0,T],\ea\ee
and
\bel{BSVIE2*}\ba{ll}
\ns\ds Y(t)=\psi(t,X(t),X(T))+\int_t^Tg(t,s,X(t),X(s),Y(s),Z(t,s),Z(s,t))ds-\int_t^T
Z(t,s)dW(s),\\
\ns\ds\qq\qq\qq\qq\qq\qq\qq\qq\qq\qq\qq\qq\qq\qq\qq t\in[0,T],\ea\ee
with $X(\cd)$ being the solution to the FSDE (\ref{SDE2}), and $\psi,g$ being some deterministic maps. Note that (\ref{BSVIE1**}) and (\ref{BSVIE2*}) are respectively Type-I and Type-II BSVIEs with random coefficients, for which the randomness all comes from the solution $X(\cd)$ of FSDE (\ref{SDE2}). Our goal is to establish a representation of $(Y(\cd),Z(\cd\,,\cd))$ in terms of $X(\cd)$, via the solution to a suitable representation PDE system. More
precisely, we will establish the following result.

\bt{Representation} \sl Let some suitable conditions hold.

\ms

{\rm(i)} For Type-I BSVIE \rf{BSVIE1**}, the following representation holds:
\bel{Rep0}\left\{\2n\ba{ll}
\ds Y(s)=\Th(s,s,X(s),X(s)),\qq\qq s\in[0,T],\\
\ns\ds Z(t,s)=\Th_x(t,s,X(t),X(s))\si(s,X(s)),\qq(t,s)\in\D[0,T],\ea\right.\ee
with $\Th(\cd\,,\cd\,,\cd\,,\cd)$ being the solution to the following PDE system:
\bel{PDE0}\left\{\2n\ba{ll}
\ds\Th_s(t,s,\xi,x)+{1\over2}\,\si(s,x)^\top\Th_{xx}(t,s,\xi,x)\si(s,x)
+\Th_x(t,s,\xi,x)b(s,x)\\
\ns\ds\qq\qq+g(t,s,\xi,x,\Th(s,s,x,x),\Th_x(t,s,\xi,x)\si(s,x)\big)=0,\qq(t,s,\xi,x)
\in\D[0,T]\times\dbR^n\times\dbR^n,\\
\ns\ds\Th(t,T,\xi,x)=\psi(t,\xi,x),\qq(t,\xi,x)\in[0,T]\times\dbR^n\times\dbR^n.\ea\right.\ee

{\rm(ii)} For Type-II BSVIE \rf{BSVIE2*}, the following representation holds:
\bel{Rep**-M}\left\{\2n\ba{ll}
\ds Y(s)=\Th(s,s,X(s),X(s)),\qq s\in[0,T],\\ [2mm]
\ns\ds Z(t,s)=\Th_x(t,s,X(t),X(s))\si(s,X(s)),\qq0\les t\les s\les T,\\ [2mm]
\ns\ds Z(t,s)=\G_x(t,s,X(s))\si(s,X(s)),\qq0\les s\les t\les T,\ea\right.\ee
with $(\G,\Th)$ being the solution to the following PDE system:
\bel{PDE**-M}\left\{\2n\ba{ll}
\ds\G_s(t,s,x)+{1\over2}\,\si(s,x)^\top\G_{xx}(t,s,x)\si(s,x)+\G_x(t,s,x)b(s,x)=0,\q0\les s\les t\les T,~x\in\dbR^n,\\
\ns\ds\Th_s(t,s,\xi,x)+{1\over2}\,\si(s,x)^\top\Th_{xx}(t,s,\xi,x)\si(s,x)+\Th_x(t,s,\xi,x)b(s,x)\\
\ns\ds\qq+g(t,s,\xi,x,\Th(s,s,x,x),\Th_x(t,s,\xi,x)\si(s,x),\G_x(s,t,\xi)\si(s,x)\big)=0,\\
\ns\ds\qq\qq\qq\qq\qq\qq\qq\qq0\les t\les s\les T,\q x,\xi\in\dbR^n,\\
\ns\ds\G(t,t,x)=\Th(t,t,x,x),\qq(t,x)\in[0,T]\times\dbR^n,\\
\ns\ds\Th(t,T,\xi,x)=\psi(t,\xi,x),\qq(t,\xi,x)\in[0,T]\times\dbR^n\times\dbR^n.\ea\right.\ee

\et

\ms

The idea of obtaining the representation is to find a proper approximation of the BSVIE by BSDEs and then derive the correct form of the representation PDE system by the invariant embedding/decoupling technique. Once the correct form of PDE system is obtained, a standard application of It\^o's formula will lead to our representation. Partial results for Type-II BSVIEs of this paper was announced in \cite{Yong 2016} without detailed proofs.

\ms

The rest of the paper is organized as follows. In Section 2, approximation of Type-I BSVIEs by means of BSDEs is established. Section 3 is devoted to the derivation of representation for the adapted solutions of Type-I BSVIEs. In Section 4, we establish a representation for the adapted M-solutions of Type-II BSVIEs. The well-posedness of representation PDEs is established in Section 5. Some concluding remarks are collected in Section 6. 

\section{Approximation of Type-I BSVIEs}

This section is devoted to an approximation of Type-I BSVIEs by a sequence of BSDEs. On one hand, such an approximation will be helpful for us to derive the representation of the adapted solutions to the Type-I BSVIEs. On the other hand, this will also be helpful for designing numerical scheme for such kind of BSVIEs (\cite{Wang 2016}). Before going further, let us first introduce the following assumption concerning the FSDE (\ref{SDE2}).

\ms

{\bf(H1)} The maps $b,\si:[0,T]\times\dbR^n\to\dbR^n$ are continuous such that for some constant $L>0$, it holds
$$|b(t,x)-b(t,x')|+|\si(t,x)-\si(t,x')|\les L|x-x'|,\qq\forall t\in[0,T],~x,x'\in\dbR^n.$$

\ms

It is standard that for any fixed $x\in\dbR^n$, FSDE (\ref{SDE2}) admits a unique strong solution $X(\cd)\equiv X(\cd\,;t,x)$, and the following holds:
\bel{|X-X|<K_1}\dbE|X(s)-X(t)|^2\les K_0|s-t|,\qq\forall s,t\in[0,T],\ee
for some constant $K_0>0$ independent of $s$ and $t$. Now, for such a given $X(\cd)$, we consider the following Type-I BSVIE:
\bel{BSVIE3}Y(t)=\psi(t,X(t),X(T))+\int_t^Tg(t,s,X(t),X(s),Y(s),Z(t,s))ds-\int_t^TZ(t,s)dW(s),\q t\in[0,T].\ee
We introduce the following assumption, recalling the definition of $\D[0,T]$ (see (\ref{D})).

\ms

{\bf(H2)} Functions $\psi:[0,T]\times\dbR^n\times\dbR^n\to\dbR^m$ and $g:\D[0,T]\times\dbR^n\times\dbR^n\times\dbR^m\times\dbR^{m\times d}\to\dbR^m$ are continuous.
There exists a constant $L>0$ such that
$$\ba{ll}
\ns\ds|\psi(t,\xi,x)-\psi(t',\xi',x')|\les L\big(|t-t'|+|\xi-\xi'|+|x-x'|\big),\q\forall(t,\xi,x),(t',\xi',x')\in[0,T]\times\dbR^n\times\dbR^n,\\
\ns\ds|g(t,s,\xi,x,y,z)-g(t',s,\xi',x',y',z')|\les L\big(|t-t'|+|\xi-\xi'|+|x-x'|+|y-y'|+|z-z'|\big),\\
\ns\ds\qq\qq\qq\qq\forall(t,s,\xi,x,y,z),(t',s,\xi',x',y',z')\1n\in\1n\D[0,T]\times
\dbR^n\times\dbR^n\times\dbR^m\times\dbR^{m\times d}.\ea$$

\ms

Under (H2), for the given $X(\cd)\equiv X(\cd\,;0,x)$, BSVIE (\ref{BSVIE3}) admits a unique adapted solution $(Y(\cd),Z(\cd\,,\cd))$ on $[0,T]$ (see \cite{Yong 2008}, for example). Let $\cP[0,T]$ be the set of all partitions $\Pi$ of $[0,T]$ having the following form:
\bel{partition}\Pi:\q0=t_0<t_1<t_2<\cds<t_{N-1}<t_N=T,\ee
with some natural number $N>1$. We define the {\it mesh size} $\|\Pi\|$ of $\Pi$ by the following:
$$\|\Pi\|=\max_{1\les i\les N}(t_i-t_{i-1}).$$
For a partition $\Pi$ as above, let us make an observation. Keep in mind that when we discuss Type-I BSVIE (\ref{BSVIE3}), the process $X(\cd)$ is given. Suppose $Y(s)$ and $Z(t,s)$ have been determined for $t_{k+1}\les t\les s\les T$ (see the region marked $\textcircled{1}$ in the figure below). Then for $t\in[t_k,t_{k+1})$, one has
$$\ba{ll}
\ns\ds Y(t)=\psi(t,X(t),X(T))+\int_t^Tg(t,s,X(t),X(s),Y(s),Z(t,s))ds
-\int_t^TZ(t,s)dW(s)\\
\ns\ds\qq=\psi(t,X(t),X(T))+\int_{t_{k+1}}^Tg(t,s,X(t),X(s),Y(s),Z(t,s))ds\\
\ns\ds\qq\qq+\int_t^{t_{k+1}}g(t,s,X(t),X(s),Y(s),Z(t,s))ds-\int_t^TZ(t,s)dW(s).\ea$$
In the above, we see that $Y(s)$, $s\in[t_{k+1},T]$, is known. But, at the moment, $Z(t,s)$ has been determined only for $t_{k+1}\les t\les s\les T$, and it is still unknown for $(t,s)\in[t_k,t_{k+1})\times[t_{k+1},T]$ (see the region marked $\textcircled{2}$ in the figure below). Also, both $Y(s)$ and $Z(t,s)$ are unknown for $t_k\les t\les s\les t_{k+1}$ (the region marked $\textcircled{3}$ in the figure below). Therefore, we want to find $(Y(s),Z(t,s))$ for $t_k\les t\les t_{k+1}$ and $t\les s\les T$ (the regions marked $\textcircled{2}$ and $\textcircled{3}$).

\vskip-0.5cm

\setlength{\unitlength}{.01in}
~~~~~~~~~~~~~~~~~~~~~~~~~~~~~~~~~~~~~~~~~~\begin{picture}(290,270)
\put(0,0){\vector(1,0){240}}
\put(0,0){\vector(0,1){240}}
\put(50,0){\line(0,1){200}}
\put(100,0){\line(0,1){200}}
\put(150,0){\line(0,1){200}}
\put(200,0){\line(0,1){200}}
\put(0,50){\line(1,0){200}}
\put(0,100){\line(1,0){200}}
\put(0,150){\line(1,0){200}}
\put(0,200){\line(1,0){200}}
\thicklines
\put(0,0){\color{red}\line(1,1){200}}
\put(-20,50){\makebox(0,0)[b]{$t_{k-1}$}}
\put(200,240){\makebox(0,0){$Y(s)$ known}}
\put(200,225){\makebox(0,0){$Z(t,s)$ known}}
\put(172,185){\makebox(0,0){$\textcircled{1}$}}
\put(175,215){\vector(0,-1){25}}
\put(100,240){\makebox(0,0){$Y(s)$ known}}
\put(100,225){\makebox(0,0){$Z(t,s)$ unknown}}
\put(128,168){\makebox(0,0){$\textcircled{2}$}}
\put(120,215){\vector(0,-1){35}}
\put(295,130){\makebox(0,0){$Y(s)$ unknown}}
\put(295,110){\makebox(0,0){$Z(t,s)$ unknown}}
\put(116,136){\makebox(0,0){$\textcircled{3}$}}
\put(240,125){\vector(-1,0){130}}
\put(-10,100){\makebox(0,0)[b]{$t_k$}}
\put(-20,150){\makebox(0,0)[b]{$t_{k+1}$}}
\put(-10,230){\makebox(0,0)[b]{$T$}}
\put(50,-22){\makebox(0,0)[b]{$t_{k-1}$}}
\put(100,-22){\makebox(0,0)[b]{$t_k$}}
\put(150,-22){\makebox(0,0)[b]{$t_{k+1}$}}
\put(230,-22){\makebox(0,0)[b]{$T$}}
\put(250,-5){\makebox{$t$}}
\put(0,245){\makebox{$s$}}
\end{picture}

\bs

\bs

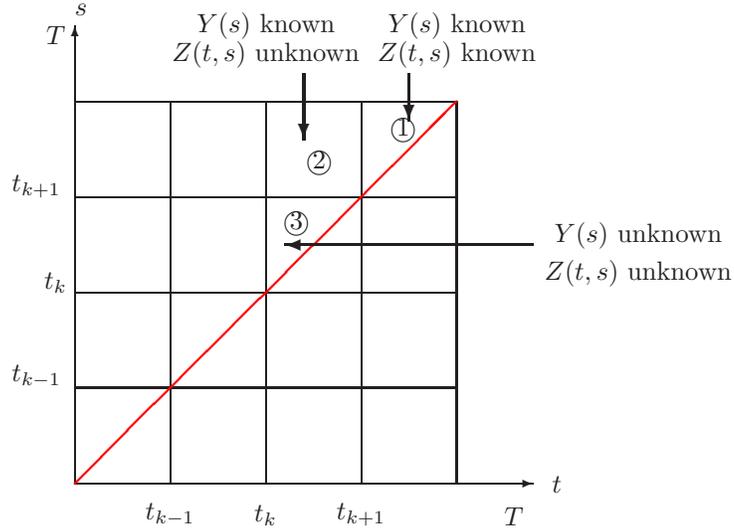
\centerline{(Figure 1)}

\bs

We now construct an approximation of BSVIE (\ref{BSVIE3}). On $[t_{N-1},T]$, we introduce the following BSDE:
\bel{BSDE(N-1)}\left\{\2n\ba{ll}
\ns\ds dY^{N-1}(s)\1n=-g\big(t_{N-1},s,X(t_{N-1}),X(s),Y^{N-1}(s),Z^{N-1}(s)\big)ds+Z^{N-1}(s)dW(s),\\
\ns\ds\qq\qq\qq\qq\qq\qq\qq\qq\qq\qq\qq\qq\qq s\in[t_{N-1},T),\\
\ns\ds Y^{N-1}(T)=\psi\big(t_{N-1},X(t_{N-1}),X(T)\big).\ea\right.\ee
Under (H2), the above BSDE admits a unique adapted solution $(Y^{N-1}(\cd),Z^{N-1}(\cd))$.
Next, on $[t_{N-2},T]$, (not just on $[t_{N-2}$, $t_{N-1}]$), we introduce the following BSDE:
\bel{BSDE(N-2)}\left\{\2n\ba{ll}
\ns\ds dY^{N-2}(s)=-g\big(t_{N-2},s,X(t_{N-2}),X(s),Y^{N-1}(s),Z^{N-2}(s)
\big)ds+Z^{N-2}(s)dW(s),\\
\ns\ds\qq\qq\qq\qq\qq\qq\qq\qq\qq\qq\qq\qq\qq s\in[t_{N-1},T),\\
\ns\ds dY^{N-2}(s)=-g\big(t_{N-2},s,X(t_{N-2}),X(s),Y^{N-2}(s),Z^{N-2}(s)
\big)ds+Z^{N-2}(s)dW(s),\\
\ns\ds\qq\qq\qq\qq\qq\qq\qq\qq\qq\qq\qq\qq\qq s\1n\in\1n[t_{N-2},t_{N-1}),\\
\ns\ds Y^{N-2}(T)=\psi\big(t_{N-2},X(t_{N-2}),X(T)\big),\qq Y^{N-2}(t_{N-1})=Y^{N-2}(t_{N-1}+0).\ea\right.\ee
Note that in the above, $Y^{N-1}(s)$ is already determined by (\ref {BSDE(N-1)}) on $[t_{N-1},T]$, which has to stay unchanged. But, since $t_{N-2}$ appears in $g$ and $\psi$, $(Y^{N-2}(\cd),Z^{N-2}(\cd))$ and $(Y^{N-1}(\cd),Z^{N-1}(\cd))$ satisfy different BSDEs on $[t_{N-1},T]$. Consequently, $Z^{N-2}(\cd)$ is possibly different from $Z^{N-1}(\cd)$ on $[t_{N-1},T]$, and therefore, $Y^{N-2}(\cd)$ could be different from $Y^{N-1}(\cd)$ on $[t_{N-1},T]$. Under our conditions, the above BSDE admits a unique adapted solution on $[t_{N-2},T]$.

\ms

In general, on $[t_k,T]$, we have a unique adapted solution $(Y^k(\cd),Z^k(\cd))$ to the  following BSDE:
\bel{BSDE(k)}\left\{\2n\ba{ll}
\ns\ds dY^k(s)=-g\big(t_k,s,X(t_k),X(s),Y^\ell(s),Z^k(s)\big)ds+Z^k(s)dW(s),\\
\ns\ds\qq\qq\qq\qq\qq\qq\qq\qq s\in[t_\ell,t_{\ell+1}),\q
k+1\les\ell\les N-1,\\
\ns\ds dY^k(s)=-g\big(t_k,s,X(t_k),X(s),Y^k(s),Z^k(s)\big)ds+Z^k(s)dW(s),\qq s\in[t_k,t_{k+1}),\\
\ns\ds Y^k(T)=\psi\big(t_k,X(t_k),X(T)\big),\q Y^k(t_\ell)=Y^k(t_\ell+0),\q k+1\les\ell\les N-1.\ea\right.\ee
Define
$$\left\{\2n\ba{ll}
\ns\ds Y^\Pi(s)=\sum_{k=0}^{N-2}Y^k(s)I_{[t_k,t_{k+1})}(s)+Y^{N-1}(s)I_{[t_{N-1},T]}(s),\qq s\in[0,T],\\
\ns\ds Z^\Pi(t,s)=\sum_{k=0}^{N-2}Z^k(s)I_{[t_k,t_{k+1})}(t)+Z^{N-1}(s)I_{[t_{N-1},T]}(t),\qq
0\les t\les s\les T,\ea\right.$$
The above defined $Y^\Pi(\cd)$ possibly has discontinuity at
$t_{N-1},t_{N-2},\cds,t_1$. With the above definition, we may rewrite (\ref{BSDE(k)}) as
\bel{BSDE(k)*}\left\{\2n\ba{ll}
\ns\ds dY^k(s)=-g\big(t_k,s,X(t_k),X(s),Y^\Pi(s),Z^k(s)\big)ds+Z^k(s)dW(s),
\qq s\in[t_k,T],\\
\ns\ds Y^k(T)=\psi\big(t_k,X(t_k),X(T)\big),\q Y^k(t_\ell)=Y^k(t_\ell+0),\q k+1\les\ell\les N-1.\ea\right.\ee
Keep in mind that for each $k=N,N-1,N-2,\cds,2,1$, process $Y^k(\cd)$ is continuous, although $Y^\Pi(\cd)$ is not necessarily continuous. Now, we introduce
\bel{t}\t^\Pi(t)=\sum_{k=0}^{N-2}t_kI_{[t_k,t_{k+1})}(t)+t_{N-1}I_{[t_{N-1},T]}(t),\qq t\in[0,T].\ee
Then
\bel{0<t-tau}0\les t-\t^\Pi(t)\les\|\Pi\|,\qq t\in[0,T].\ee
Hence, for any $t\in[0,T)$, let $t\in[t_k,t_{k+1})$, we have the following:
$$\ba{ll}
\ns\ds Y^\Pi(t)=Y^k(t)=\psi\big(t_k,X(t_k),X(T)\big)+\int_t^{t_{k+1}}g\big(t_k,
s,X(t_k),X(s),Y^k(s),Z^k(s)\big)ds\\
\ns\ds\qq\qq\qq\qq+\sum_{\ell=k+1}^{N-1}\int_{t_\ell}^{t_{\ell+1}}g\big(t_k,s,X(t_k),
X(s),Y^\ell(s),Z^k(s)\big)ds-\int_t^TZ^k(s)dW(s)\\
\ns\ds\qq=\psi\big(\t^\Pi(t),X(\t^\Pi(t)),X(T)\big)+\int_t^Tg\big(\t^\Pi(t),s,
X(\t^\Pi(t)),X(s),Y^\Pi(s),Z^\Pi(t,s)\big)ds-\int_t^TZ^\Pi(t,s)dW(s).\ea$$
Therefore, $(Y^\Pi(\cd),Z^\Pi(\cd\,,\cd))$ satisfies a BSDE on each $(t_k,t_{k+1})$, and satisfies a BSVIE on $[0,T]$. Note that since $t\mapsto g(\t^\Pi(t),s,X(t),x,y,z)$ has possible jumps at $t_k$, the resulting $t\mapsto Y^\Pi(t)$ may also have jumps at $t_k$, regardless of its integral form. For the above constructed $(Y^\Pi(\cd),Z^\Pi(\cd\,,\cd))$, we have the following proposition.

\bp{Proposition 2.1.} \sl Suppose {\rm(H1)--(H2)} hold, $X(\cd)$ is the solution to the FSDE $(\ref{SDE2})$, and $(Y^\Pi(\cd),Z^\Pi(\cd\,,\cd))$ is constructed as above. Then
\bel{Y-Y}\dbE\int_0^T|Y^\Pi(t)-Y(t)|^2dt+\dbE\int_0^T\3n\int_t^T|Z^\Pi(t,s)-Z(t,s)|^2dsdt\les
K\|\Pi\|^2,\ee
and
\bel{Y(k+1)-Y(k)}\dbE\[\sup_{s\in[t_{k+1},T]}|Y^{k+1}(s)-Y^k(s)|^2\]
+\dbE\(\int_{t_{k+1}}^T|Z^{k+1}(s)-Z^k(s)|^2ds\)\les K\|\Pi\|,\q0\les k\les N-2.\ee
In particular,
\bel{}\dbE|Y^\Pi(t_k)-Y^\Pi(t_k-0)|\les K\|\Pi\|,\qq1\les k\les N-1.\ee

\ep

\it Proof. \rm By the stability of adapted solutions to BSVIEs (\cite{Yong 2008}), we have (note \eqref{0<t-tau})
$$\ba{ll}
\ns\ds\dbE\int_0^T|Y^\Pi(t)-Y(t)|^2dt+\dbE\int_0^T\3n\int_t^T|Z^\Pi(t,s)-Z(t,s)|^2dsdt\\
\ns\ds\les K\dbE\int_0^T|\psi\big(\t^\Pi(t),X(t),X(T)\big)-\psi(t,X(t),X(T))|^2dt\\
\ns\ds\qq+\dbE\int_0^T\3n\int_t^T|g\big(\t^\Pi(t),s,X(t),X(s),Y(s),Z(t,s)\big)
-g\big(t,s,X(t),X(s),Y(s),Z(t,s)\big)|^2dsdt\\
\ns\ds\les K\dbE\Big\{\int_0^T\|\Pi\|^2dt+\int_0^T\3n\int_t^T\|\Pi\|^2dsdt
\Big\}\les K\|\Pi\|^2.\ea$$
This proves (\ref{Y-Y}). Next, one has
$$\left\{\2n\ba{ll}
\ns\ds d\[Y^{k+1}(s)-Y^k(s)\]=-\[g(t_{k+1},s,X(t_{k+1}),X(s),Y^{k+1}(s),Z^{k+1}(s))\\
\ns\ds\qq\qq\qq\qq\qq\qq-g(t_k,s,X(t_k),X(s),Y^k(s),Z^k(s))\]ds
+\[Z^{k+1}(s)-Z^k(s)\]dW(s),\\
\ns\ds\qq\qq\qq\qq\qq\qq\qq\qq\qq\qq\qq\qq\qq\qq\qq s\in[t_{k+1},T],\\
\ns\ds\[Y^{k+1}(T)-Y^k(T)\]=\psi(t_{k+1},X(t_{k+1}),X(T))-\psi(t_k,X(t_k),X(T)).
\ea\right.$$
Hence,
$$\ba{ll}
\ns\ds\dbE\[\sup_{s\in[t_{k+1},T]}|Y^{k+1}(s)-Y^k(s)|^2\]
+\dbE\(\int_{t_{k+1}}^T|Z^{k+1}(s)-Z^k(s)|^2ds\)\\
\ns\ds\les K\dbE\Big\{|\psi(t_{k+1},X(t_{k+1}),X(T))-\psi(t_k,X(t_k),X(T))|^2\\
\ns\ds\q+\(\int_{t_{k+1}}^T|g(t_{k+1},s,X(t_{k+1}),X(s),Y^k(s),Z^k(s))
-g(t_k,s,X(t_k),X(s),Y^k(s),Z^k(s))|ds\)^2\Big\}\\
\ns\ds\les K\dbE\Big\{(t_{k+1}-t_k)^2+|X(t_{k+1})-X(t_k)|^2+\int_{t_{k+1}}^T\(|t_{k+1}-t_k|
+|X(t_{k+1})-X(t_k)|\)ds\Big\}^2\les K\|\Pi\|.\ea$$
This leads to our conclusion. \endpf

\ms

\section{Representation of Adapted Solutions for Type-I BSVIEs}

In this section, we will represent the adapted solution $(Y(\cd),Z(\cd\,,\cd))$ of Type-I BSVIE (\ref{BSVIE3}) in terms of $X(\cd)$, the solution to FSDE (\ref{SDE2}), and the solution to the corresponding representation PDE.

\ms

Let $X(\cd)$ be the solution of \eqref{SDE2} and $\Pi\in\cP[0,T]$ be of form (\ref{partition}). Let  $(Y^k(\cd),Z^k(\cd))$ ($0\les k\les N-1$) and $(Y^\Pi(\cd),Z^\Pi(\cd\,,\cd))$ be constructed as in Section 2.

\ms

If $(Y^\Pi(\cd),Z^\Pi(\cd\,,\cd))$ were represented by $X(\cd)$, together with the solution to certain PDE, then by sending $\|\Pi\|\to0$, we would get what we want. However, there are some difficulties in doing that directly (see below for some explanations). Therefore, instead, we construct a sequence of processes $(\bar Y^\Pi(\cd),\bar Z^\Pi(\cd\,,\cd))$ which is close to $(Y^\Pi(\cd),Z^\Pi(\cd\,,\cd))$ and which can be represented by $X(\cd)$, together with the solution of a certain PDE. Then by sending $\|\Pi\|\to0$, we will obtain the desired representation of $(Y(\cd),Z(\cd\,,\cd))$.

\ms

Now, we carefully make this precise. First, let $(\bar Y^{N-1}(\cd),\bar Z^{N-1}(\cd))$ be the adapted solution to the following BSDE:
\bel{FBSDE(N-1)}\left\{\2n\ba{ll}
\ns\ds d\bar Y^{N-1}(s)=-g\big(t_{N-1},s,X(t_{N-1}),X(s),\bar Y^{N-1}(s),\bar Z^{N-1}(s)\big)ds\1n+\1n\bar Z^{N-1}(s)dW(s),\qq s\in[t_{N-1},T],\\
\ns\ds\bar Y^{N-1}(T)=\psi\big(t_{N-1},X(t_{N-1}),X(T)\big),\ea\right.\ee
which coincides with BSDE (\ref{BSDE(N-1)}). Thus, one has
\bel{YZ=YZ}(\bar Y^{N-1}(s),\bar Z^{N-1}(s))=(Y^{N-1}(s),Z^{N-1}(s)),\qq t_{N-1}\les t\les s\les T.\ee
Although in (\ref{FBSDE(N-1)}), the map $(s,x,y,z)\mapsto\big(g(t_{N-1},s,X(t_{N-1}),x,y,z),
\psi(t_{N-1},X(t_{N-1}),x)\big)$ is merely $\cF_{t_{N-1}}$-measurable, not necessarily deterministic, since we are considering the BSDE on $[t_{N-1},T]$, the (decoupling) technique introduced in \cite{Ma-Protter-Yong 1994, Ma-Yong 1999} will still work. In fact, we have the following representation:
\bel{Rep-YZ(N-1)}\left\{\2n\ba{ll}
\ns\ds\bar Y^{N-1}(s)=Y^{N-1}(s)=\Th^{N-1}\big(s,X(t_{N-1}),X(s)\big),\\
\ns\ds\bar Z^{N-1}(s)=Z^{N-1}(s)=\Th^{N-1}_x\big(s,X(t_{N-1}),X(s)\big)\si(s,X(s)),\ea\right.\qq s\in[t_{N-1},T],\ee
with $(s,x)\mapsto\Th^{N-1}(s,\xi,x)$ being the solution to the following representation PDE:
\bel{Th(N-1)}\left\{\2n\ba{ll}
\ns\ds\Th^{N-1}_s(s,\xi,x)+{1\over2}\,\si(s,x)^\top\Th^{N-1}_{xx}(s,\xi,x)
\si(s,x)+\Th^{N-1}_x(s,\xi,x)b(s,x)\\
\ns\ds\qq\q+g\big(t_{N-1},s,\xi,x,\Th^{N-1}(s,\xi,x),\Th^{N-1}_x(s,\xi,x)
\si(s,x)\big)=0,\qq
(s,x)\in[t_{N-1},T]\times\dbR^n,\\
\ns\ds\Th^{N-1}(T,\xi,x)=\psi(t_{N-1},\xi,x),\qq x\in\dbR^n.\ea\right.\ee
In the above, $\xi\in\dbR^n$ is treated as a parameter. With the representation \eqref{Rep-YZ(N-1)}, we can rewrite (\ref{FBSDE(N-1)}) as follows:
\bel{BSDE(N-1)*}\left\{\2n\ba{ll}
\ns\ds d\bar Y^{N-1}(s)=-g\big(t_{N-1},s,X(t_{N-1}),X(s),\Th^{N-1}(s,X(t_{N-1}),X(s)),
\bar Z^{N-1}(s)\big)ds\\
\ns\ds\qq\qq\qq\qq\qq\qq\qq\qq\qq\qq\qq+\bar Z^{N-1}(s)dW(s),\qq s\in[t_{N-1},T),\\
\ns\ds\bar Y^{N-1}(T)=\psi\big(t_{N-1},X(t_{N-1}),X(T)\big).\ea\right.\ee

\ms

Next, we construct $(\bar Y^{N-2}(\cd),\bar Z^{N-2}(\cd))$ on $[t_{N-2},T]$. On $[t_{N-1},T]$, we let $(\bar Y^{N-2}(\cd),\bar Z^{N-2}(\cd))$ be the adapted solution to the following BSDE:
\bel{FBSDE(N-2)}\left\{\2n\ba{ll}
\ns\ds d\bar Y^{N-2}(s)=-g\big(t_{N-2},s,X(t_{N-2}),X(s),\Th^{N-1}(s,X(s),X(s)),\bar Z^{N-2}(s)\big)ds\\
\ns\ds\qq\qq\qq\qq\qq\qq\qq\qq\qq\qq+\bar Z^{N-2}(s)dW(s),\qq s\in[t_{N-1},T],\\
\ns\ds\bar Y^{N-2}(T)=\psi\big(t_{N-2},X(t_{N-2}),X(T)\big).\ea\right.\ee
Note that on $[t_{N-1},T]$, we have \rf{YZ=YZ} and representation \rf{Rep-YZ(N-1)}.
Hence, by \rf{BSDE(N-2)}, we see that $(Y^{N-2}(\cd),Z^{N-2}(\cd))$ satisfies the following BSDE:
\bel{BSDE(N-2)*}\left\{\2n\ba{ll}
\ns\ds dY^{N-2}(s)=-g\big(t_{N-2},s,X(t_{N-2}),X(s),\Th^{N-1}(s,X(t_{N-1}),X(s)),
Z^{N-2}(s)\big)ds\\
\ns\ds\qq\qq\qq\qq\qq\qq\qq\qq\qq\qq\qq+Z^{N-2}(s)dW(s),\qq s\in[t_{N-1},T],\\
\ns\ds Y^{N-2}(T)=\psi\big(t_{N-2},X(t_{N-2}),X(T)\big).\ea\right.\ee
Let us make two comparisons. First, \eqref{BSDE(N-1)*} and \eqref{BSDE(N-2)*} are different: $(t_{N-1},X(t_{N-1}))$ in the former is replaced by $(t_{N-2},X(t_{N-2}))$ in the latter at two places. Second, \eqref{BSDE(N-2)*} and \eqref{FBSDE(N-2)} are different: $\Th^{N-1}(s,X(t_{N-1}),$ $X(s))$ in the former is replaced by $\Th^{N-1}(s,X(s),X(s))$ in the latter. Note that in (\ref{BSDE(N-2)*}), both $X(t_{N-1})$ and $X(t_{N-2})$ appear. This will cause some difficulties in passing to the limit as $\|\Pi\|\to0$ later on. This is exactly the difficulty that we will encounter if we use $(Y^\Pi(\cd),Z^\Pi(\cd))$ directly trying to get our representation. On the other hand, since $\|\Pi\|$ will be small, $X(t_{N-1})$ and $X(s)$ will be close (in some sense), for $s\in[t_{N-1},T]$, it should be harmless to replace $\Th^{N-1}\big(s,X(t_{N-1}),X(s)\big)$ by $\Th^{N-1}\big(s,X(s),X(s)\big)$ in the drift of the equation on $[t_{N-1},T]$.
Thus, if $\|\Th^{N-1}_{\xi}\|_{\infty}^2<\infty$, by the stability of adapted solutions to BSDEs, we have
\bel{4.8}\ba{ll}
\ns\ds\dbE\[\sup_{s\in[t_{N-1},T]}|\bar Y^{N-2}(s)-Y^{N-2}(s)|^2\]
+\dbE\int_{t_{N-1}}^T|\bar Z^{N-2}(s)-Z^{N-2}(s)|^2ds\\
%
%
\ns\ds\les K_1\dbE\(\int_{t_{N-1}}^T\big|g(t_{N-2},s,X(t_{N-2}),X(s),\Th^{N-1}(s,X(t_{N-1}),
X(s)),Z^{N-2}(s))\\
\ns\ds\qq\qq\qq-g(t_{N-2},s,X(t_{N-2}),X(s),\Th^{N-1}(s,X(s),X(s)),Z^{N-2}(s))\big|ds\)^2\\
\ns\ds\les K_1L^2\dbE\(\int_{t_{N-1}}^T|\Th^{N-1}(s,X(s),X(s))
-\Th^{N-1}(s,X(t_{N-1}),X(s))|ds\)^2\\
\ns\ds\les K_1L^2\|\Th^{N-1}_\xi\|_\infty^2\dbE\(\int_{t_{N-1}}^T|X(s)
-X(t_{N-1})|ds\)^2\\
\ns\ds\les K_1L^2\|\Th^{N-1}_\xi\|_\infty^2(T-t_{N-1})\int_{t_{N-1}}^T\dbE|X(s)-X(t_{N-1})|^2ds\\
\ns\ds\les K_1L^2\|\Th^{N-1}_\xi\|_\infty^2K_0(T-t_{N-1})\int_{t_{N-1}}^T(s-t_{N-1})ds
\les{K_0K_1L^2\|\Th^{N-1}_\xi\|_\infty^2\|\Pi\|\over2}(T-t_{N-1})^2.\ea\ee
In the above, $K_0$ is the constant appears in (\ref{|X-X|<K_1}), and $K_1$ is a constant appears in the stability estimate for the adapted solution of BSDEs, which can be made independent of the partition $\Pi$, under (H2).

\ms

Similar to the previous step, for (\ref{FBSDE(N-2)}), we have the following representation:
\bel{Rep-YZ(N-2)}\left\{\2n\ba{ll}
\ns\ds\bar Y^{N-2}(s)=\Th^{N-2}\big(s,X(t_{N-2}),X(s)\big),\\
\ns\ds\bar Z^{N-2}(s)=\Th^{N-2}_x\big(s,X(t_{N-2}),X(s)\big)\si(s,X(s)),\ea\right.\qq s\in[t_{N-1},T],\ee
with $(s,x)\mapsto\Th^{N-2}(s,\xi,x)$ being the solution to the following representation PDE:
\bel{Th(N-2)}\left\{\2n\ba{ll}
\ns\ds\Th^{N-2}_s(s,\xi,x)+{1\over2}\,\si(s,x)^\top\Th^{N-2}_{xx}(s,\xi,x)
\si(s,x)+\Th^{N-2}_x(s,\xi,x)b(s,x)\\
\ns\ds\qq\q+g\big(t_{N-2},s,\xi,x,\Th^{N-1}(s,x,x),\Th^{N-2}_x(s,\xi,x)
\si(s,x)\big)=0,\qq(s,x)\in[t_{N-1},T]\times\dbR^n,\\
\ns\ds\Th^{N-2}(T,\xi,x)=\psi(t_{N-2},\xi,x),\qq x\in\dbR^n.\ea\right.\ee
Note that equation (\ref{Th(N-2)}) is different from (\ref{Th(N-1)}), not just because $t_{N-1}$ is replaced by $t_{N-2}$ in $g$ and $\psi$, but also because $\Th^{N-1}(s,\xi,x)$ is replaced by $\Th^{N-1}(s,x,x)$ in $g$. We expect that $\Th^{N-2}(s,\xi,x)\big|_{\xi=x}$ is close to $\Th^{N-1}(s,\xi,x)\big|_{\xi=x}$, for $s\in[t_{N-1},T]$, when $t_{N-1}-t_{N-2}>0$ is small.

\ms

So far, we have constructed $(\bar Y^{N-2}(\cd),\bar Z^{N-2}(\cd)$ on $[t_{N-1},T]$. To construct $(\bar Y^{N-2}(\cd),\bar Z^{N-2}(\cd))$ on $[t_{N-2},t_{N-1})$, we introduce the following BSDE on $[t_{N-2},t_{N-1})$:
\bel{FBSDE(N-2)*}\left\{\2n\ba{ll}
\ns\ds d\bar Y^{N-2}(s)=-g(t_{N-2},s,X(t_{N-2}),X(s),\bar Y^{N-2}(s),\bar Z^{N-2}(s))ds+\bar Z^{N-2}(s)dW(s),\\
\ns\ds\qq\qq\qq\qq\qq\qq\qq\qq\qq\qq\qq\qq\qq s\in[t_{N-2},t_{N-1}),\\
\ns\ds\bar Y^{N-2}(t_{N-1}-0)=\Th^{N-2}(t_{N-1},X(t_{N-2}),X(t_{N-1})).\ea\right.\ee
Now, on $[t_{N-2},t_{N-1})$, we have the following representation:
\bel{Rep-YZ(N-2)*}\left\{\2n\ba{ll}
\ns\ds\bar Y^{N-2}(s)=\Th^{N-2}\big(s,X(t_{N-2}),X(s)\big),\\
\ns\ds\bar Z^{N-2}(s)=\Th^{N-2}_x(s,X\big(t_{N-2}),X(s)\big)\si(s,X(s)),\ea\right.\qq s\in[t_{N-2},t_{N-1}),\ee
with $(s,x)\mapsto\Th^{N-2}(s,\xi,x)$ being the solution to the following representation  PDE:
\bel{Th(N-2)**}\left\{\2n\ba{ll}
\ns\ds\Th^{N-2}_s(s,\xi,x)+{1\over2}\,\si(s,x)^\top\Th^{N-2}_{xx}(s,\xi,x)
\si(s,x)+\Th^{N-2}_x(s,\xi,x)b(s,x)\\
\ns\ds\q+g\big(t_{N-2},s,\xi,x,\Th^{N-2}(s,\xi,x),\Th^{N-2}_x(s,\xi,x)
\si(s,x)\big)=0,\q(s,x)\in[t_{N-2},t_{N-1})\times\dbR^n,\\
\ns\ds\Th^{N-2}(t_{N-1},\xi,x)=\Th^{N-2}(t_{N-1}+0,\xi,x),\qq x\in\dbR^n.\ea\right.\ee
Note that unlike \eqref{Th(N-2)}, in the above, $\Th^{N-2}(s,\xi,x)$ appears instead of $\Th^{N-2}(s,x,x)$ in $g$. Next, since $(\bar Y^{N-2}(\cd),\bar Z^{N-2}(\cd))$ and $(Y^{N-2}(\cd),Z^{N-2}(\cd))$ satisfy the same equation on $[t_{N-2},t_{N-1})$ with different terminal conditions $\bar Y^{N-2}(t_{N-1})$ and $Y^{N-2}(t_{N-1})$, we must have, making use of \eqref{4.8},
\bel{3.14}\ba{ll}
\ns\ds\dbE\[\sup_{s\in[t_{N-2},t_{N-1})}|\bar Y^{N-2}(s)-Y^{N-2}(s)|^2+
\int_{t_{N-2}}^{t_{N-1}}|\bar Z^{N-2}(s)-Z^{N-2}(s)|^2ds\]\\ [3mm]
\ns\ds\les K_0\dbE\,|\bar Y^{N-2}(t_{N-1})-Y^{N-2}(t_{N-1})|^2\les{K_0^2K_1L^2\|\Th^{N-1}_\xi\|_\infty^2\|\Pi\|\over2}
(T-t_{N-1})^2.\ea\ee
To summarize the above, we have
\bel{Rep-YZ(N-2)*}\left\{\2n\ba{ll}
\ns\ds\bar Y^{N-2}(s)=\Th^{N-2}\big(s,X(t_{N-2}),X(s)\big),\\
\ns\ds\bar Z^{N-2}(s)=\Th^{N-2}_x\big(s,X(t_{N-2}),X(s)\big)\si(s,X(s)),\ea\right.\qq s\in[t_{N-2},T],\ee
with $(s,x)\mapsto\Th^{N-2}(s,\xi,x)$ being the solution to the following:
\bel{Th(N-2)**}\left\{\2n\ba{ll}
\ns\ds\Th^{N-2}_s(s,\xi,x)+{1\over2}\,\si(s,x)^\top\Th^{N-2}_{xx}(s,\xi,x)
\si(s,x)+\Th^{N-2}_x(s,\xi,x)b(s,x)\\
\ns\ds\q+g\big(t_{N-2},s,\xi,x,\Th^{N-1}(s,x,x),\Th^{N-2}_x(s,\xi,x)
\si(s,x)\big)=0,\q(s,x)\in[t_{N-1},T]\times\dbR^n,\\
\ns\ds\Th^{N-2}_s(s,\xi,x)+{1\over2}\,\si(s,x)^\top\Th^{N-2}_{xx}(s,\xi,x)
\si(s,x)+\Th^{N-2}_x(s,\xi,x)b(s,x)\\
\ns\ds\q+g\big(t_{N-2},s,\xi,x,\Th^{N-2}(s,\xi,x),\Th^{N-2}_x(s,\xi,x)
\si(s,x)\big)=0,\q(s,x)\in[t_{N-2},t_{N-1})\times\dbR^n,\\
\ns\ds\Th^{N-2}(T,\xi,x)=\psi(t_{N-2},\xi,x),\q\Th^{N-2}(t_{N-1},\xi,x)
=\Th^{N-2}(t_{N-1}+0,\xi,x),\qq x\in\dbR^n.\ea\right.\ee
Note that although $s\mapsto\bar Y^{N-2}(s)$ could be discontinuous at $s=t_{N-1}$, the function $s\mapsto\Th^{N-2}(s,\xi,x)$ is continuous. Also, we point out that in the above system \eqref{Th(N-2)**}, the equations on $[t_{N-1},T]$ and $[t_{N-2},t_{N-1})$ are different: $\Th^{N-1}(s,x,x)$ appears in $g$ for the former and $\Th^{N-2}(s,\xi,x)$ appears in $g$
for the latter.

\ms

The above discussion seems not enough to obtain an inductive statement. In particular, we need to make sure that the estimate on the error between $(\bar Y^k(\cd),\bar Z^k(\cd))$ and $(Y^k(\cd),Z^k(\cd))$ will not be unboundedly accumulated. Thus, let us construct $(\bar Y^{N-3}(\cd),\bar Z^{N-3}(\cd))$ on $[t_{N-3},T]$. To this end, we consider the following BSDE on $[t_{N-2},T]$:
\bel{FBSDE(N-3)}\left\{\2n\ba{ll}
\ns\ds d\bar Y^{N-3}(s)=-g\big(t_{N-3},s,X(t_{N-3}),X(s),\Th^{N-1}(s,X(s),X(s)),
\bar Z^{N-3}(s)\big)ds\\
\ns\ds\qq\qq\qq\qq\qq\qq\qq\qq\qq\qq\qq+\bar Z^{N-3}(s)dW(s),\q s\in[t_{N-1},T],\\
\ns\ds d\bar Y^{N-3}(s)=-g\big(t_{N-3},s,X(t_{N-3}),X(s),\Th^{N-2}(s,X(s),X(s)),
\bar Z^{N-3}(s)\big)ds\\
\ns\ds\qq\qq\qq\qq\qq\qq\qq\qq\qq\qq\qq+\bar Z^{N-3}(s)dW(s),\q s\in[t_{N-2},t_{N-1}),\\
\ns\ds\bar Y^{N-3}(T)=\psi\big(t_{N-3},X(t_{N-3}),X(T)\big),\\
\ns\ds\bar Y^{N-3}(t_{N-1}-0)=\Th^{N-3}\big(t_{N-1},X(t_{N-3}),X(t_{N-1})\big),\ea\right.\ee
where $\Th^{N-3}(\cd\,,\cd\,,\cd)$ is the solution to the following PDE:
\bel{Th(N-3)}\left\{\2n\ba{ll}
\ns\ds\Th^{N-3}_s(s,\xi,x)+{1\over2}\,\si(s,x)^T\Th^{N-3}_{xx}(s,\xi,x)
\si(s,x)+\Th^{N-3}_x(s,\xi,x)b(s,x)\\
\ns\ds\q+g\big(t_{N-3},s,\xi,x,\Th^{N-1}(s,x,x),\Th^{N-3}_x(s,\xi,x)
\si(s,x)\big)=0,\q(s,x)\in[t_{N-1},T]\times\dbR^n,\\
\ns\ds\Th^{N-3}_s(s,\xi,x)+{1\over2}\,\si(s,x)^T\Th^{N-3}_{xx}(s,\xi,x)
\si(s,x)+\Th^{N-3}_x(s,\xi,x)b(s,x)\\
\ns\ds\q+g\big(t_{N-3},s,\xi,x,\Th^{N-2}(s,x,x),\Th^{N-3}_x(s,\xi,x)
\si(s,x)\big)=0,\q(s,x)\in[t_{N-2},t_{N-1})\times\dbR^n,\\
\ns\ds\Th^{N-3}(T,\xi,x)=\psi(t_{N-3},\xi,x),\q\Th^{N-3}(t_{N-1},\xi,x)
=\Th^{N-3}(t_{N-1}+0,\xi,x),\qq x\in\dbR^n.\ea\right.\ee
Then we have the following representation:
\bel{}\left\{\2n\ba{ll}
\ns\ds\bar Y^{N-3}(s)=\Th^{N-3}(s,X(t_{N-3}),X(s)),\\
\ns\ds\bar Z^{N-3}(s)=\Th^{N-3}_x(s,X(t_{N-3}),X(s))\si(s,X(s)),\ea\right.\qq
s\in[t_{N-2},T],\ee
By \rf{Rep-YZ(N-1)} and \rf{Rep-YZ(N-2)*}, we see that
$$\left\{\2n\ba{ll}
\ns\ds dY^{N-3}(s)=-g\big(t_{N-3},s,X(t_{N-3}),X(s),\Th^{N-1}(s,X(t_{N-1}),X(s)),
Z^{N-3}(s)\big)ds\\
\ns\ds\qq\qq\qq\qq\qq\qq\qq\qq\qq\qq\qq+Z^{N-3}(s)dW(s),\qq s\in[t_{N-1},T],\\
\ns\ds dY^{N-3}(s)=-g\big(t_{N-3},s,X(t_{N-3}),X(s),\Th^{N-2}(s,X(t_{N-2}),X(s)),
Z^{N-3}(s)\big)ds\\
\ns\ds\qq\qq\qq\qq\qq\qq\qq\qq\qq\qq\qq+Z^{N-3}(s)dW(s),\qq s\in[t_{N-2},t_{N-1}),\\
\ns\ds Y^{N-3}(T)=\psi(t_{N-3},X(t_{N-3}),X(T)).\ea\right.$$
Thus, by the stability of adapted solutions to BSDEs, one has
$$\ba{ll}
\ns\ds\dbE\[\sup_{s\in[t_{N-2},T]}|\bar Y^{N-3}(s)-Y^{N-3}(s)|^2+\int_{t_{N-2}}^T
|\bar Z^{N-3}(s)-Z^{N-3}(s)|^2ds\]\\
\ns\ds\les K_1\dbE\(\int_{t_{N-1}}^T\big|g(t_{N-3},s,X(t_{N-3}),\Th^{N-1}(s,X(s),X(s)),
Z^{N-3}(s)\big)\\
\ns\ds\qq\qq-g\big(t_{N-3},s,X(t_{N-3}),\Th^{N-1}(s,X(t_{N-1}),X(s)),
Z^{N-3}(s)\big)\big|ds\\
\ns\ds\qq\qq+\int_{t_{N-2}}^{t_{N-1}}\big|g(t_{N-3},s,X(t_{N-3}),\Th^{N-2}(s,X(s),X(s)),
Z^{N-3}(s)\big)\\
\ns\ds\qq\qq-g\big(t_{N-3},s,X(t_{N-3}),\Th^{N-2}(s,X(t_{N-2}),X(s)),
Z^{N-3}(s)\big)\big|ds\)^2\\
\ns\ds\les K_1L^2\dbE\(\int_{t_{N-1}}^T\big|\Th^{N-1}(s,X(s),X(s))-\Th^{N-1}(s,X(t_{N-1}),
X(s))\big|ds\\
\ns\ds\qq\qq+\int_{t_{N-2}}^{t_{N-1}}\big|\Th^{N-2}(s,X(s),X(s))-\Th^{N-2}(s,X(t_{N-2}),
X(s))\big|ds\)^2\\
\ns\ds\les K_1L^2\|\Th_\xi\|_\infty^2\dbE\(\int_{t_{N-1}}^T\big|X(s)-X(t_{N-1})\big|ds
+\int_{t_{N-2}}^{t_{N-1}}\big|X(s)-X(t_{N-2})\big|ds\)^2\\
\ns\ds\les K_1L^2\|\Th_\xi\|_\infty^2(T-t_{N-2})\(\int_{t_{N-1}}^T\dbE|X(s)-X(t_{N-1})|^2ds
+\int_{t_{N-2}}^{t_{N-1}}\dbE|X(s)-X(t_{N-2})|^2ds\)\\
\ns\ds\les K_1L^2\|\Th_\xi\|_\infty^2(T-t_{N-2})K_0\(\int_{t_{N-1}}^T(s-t_{N-1})ds+\int_{t_{N-2}}^{t_{N-1}}
(s-t_{N-2})ds\)\\
\ns\ds\les{K_0K_1L^2\|\Th_\xi\|_\infty^2\over2}(T-t_{N-2})\((T-t_{N-1})^2+(t_{N-1}-t_{N-2})^2\)\\
\ns\ds\les{K_0K_1L^2\|\Th_\xi\|_\infty^2\over2}(T-t_{N-2})^2\|\Pi\|,\ea$$
where
$$\|\Th_\xi\|_\infty=\|\Th_\xi^{N-1}\|_\infty\vee\|\Th_\xi^{N-2}\|_\infty.$$
Next, on $[t_{N-3},t_{N-2})$, $(\bar Y^{N-3}(\cd),\bar Z^{N-3}(\cd))$ and $(Y^{N-3}(\cd),Z^{N-3}(\cd))$ satisfy the same equation with possibly different terminal conditions at $t=t_{N-2}$. Hence, we have
\bel{3.20}\ba{ll}
\ns\ds\dbE\[\sup_{s\in[t_{N-3},t_{N-2})}|\bar Y^{N-3}(s)-Y^{N-3}(s)|^2
+\int_{t_{N-3}}^{t_{N-2}}|\bar Z^{N-3}(s)-Z^{N-3}(s)|^2ds\]\\
\ns\ds\les K_1\dbE|\bar Y^{N-3}(t_{N-2})-Y^{N-3}(t_{N-2})|^2
\les{K_0K_1^2L^2\|\Th_\xi\|_\infty^2\over2}(T-t_{N-2})^2\|\Pi\|.\ea\ee
Also, we have
$$\left\{\2n\ba{ll}
\ns\ds\bar Y^{N-3}(s)=\Th^{N-3}(s,X(t_{N-3}),X(s)),\\
\ns\ds\bar Z^{N-3}(s)=\Th^{N-3}_x(s,X(t_{N-3}),X(s))\si(s,X(s)),\ea\right.\qq s\in[t_{N-3},t_{N-2}),$$
with $\Th^{N-3}(\cd\,,\cd)$ satisfying
\bel{Th(N-3)**}\left\{\2n\ba{ll}
\ns\ds\Th^{N-3}_s(s,\xi,x)+{1\over2}\,\si(s,x)^\top\Th^{N-3}_{xx}(s,\xi,x)
\si(s,x)+\Th^{N-3}_x(s,\xi,x)b(s,x)\\
\ns\ds\q+g\big(t_{N-3},s,\xi,x,\Th^{N-3}(s,\xi,x),\Th^{N-3}_x(s,\xi,x)
\si(s,x)\big)=0,\q(s,x)\in[t_{N-3},t_{N-2})\times\dbR^n,\\
\ns\ds\Th^{N-3}(t_{N-2},\xi,x)=\Th^{N-3}(t_{N-2}+0,\xi,x),\qq x\in\dbR^n.\ea\right.\ee

\ms

Now, we look at the general case. For each $k=0,1,\cds,N-1$, on $[t_k,T]$, we consider the following BSDE:
\bel{FBSDE(k)}\left\{\2n\ba{ll}
\ns\ds d\bar Y^k(s)=-g\big(t_k,s,X(t_k),X(s),\Th^\ell(s,X(s),X(s)),\bar Z^k(s)\big)ds+\bar Z^k(s)dW(s),\\
\ns\ds\qq\qq\qq\qq\qq\qq\qq\qq\q s\in[t_\ell,t_{\ell+1}),\q k+1\les\ell\les N-1,\\
\ns\ds d\bar Y^k(s)=-g\big(t_k,s,X(t_k),X(s),\bar Y^k(s),\bar Z^k(s)\big)ds+\bar Z^k(s)dW(s),\q s\in[t_k,t_{k+1}),\\
\ns\ds\bar Y^k(T)=\psi\big(t_k,X(t_k),X(T)\big).\ea\right.\ee
Then the following representation holds:
\bel{Rep-YZ(k)}\left\{\2n\ba{ll}
\ns\ds\bar Y^k(s)=\Th^k(s,X(t_k),X(s)),\\
\ns\ds\bar Z^k(s)=\Th^k_x\big(s,X(t_k),X(s)\big)\si(s,X(s)),\ea\right.\qq s\in[t_k,T],\ee
where $(s,x)\mapsto\Th^k(s,\xi,x)$ is the solution to the following PDE:
\bel{Th(k)}\left\{\2n\ba{ll}
\ns\ds\Th^k_s(s,\xi,x)+{1\over2}\,\si(s,x)^\top\Th^k_{xx}(s,\xi,x)
\si(s,x)+\Th^k_x(s,\xi,x)b(s,x)\\
\ns\ds\qq+g\big(t_k,s,\xi,x,\Th^\ell(s,x,x),\Th^k_x(s,\xi,x)\si(s,x)
\big)=0,\\
\ns\ds\qq\qq\qq\qq\qq\qq\qq\qq\qq(s,x)\in[t_\ell,t_{\ell+1})\times\dbR^n,\q k+1\les\ell\les N-1,\\
\ns\ds\Th^k_s(s,\xi,x)+{1\over2}\,\si(s,x)^\top\Th^k_{xx}(s,\xi,x)
\si(s,x)+\Th^k_x(s,\xi,x)b(s,x)\\
\ns\ds\qq+g\big(t_k,s,\xi,x,\Th^k(s,\xi,x),\Th^k_x(s,\xi,x)\si(s,x)\big)=0,
\q(s,x)\in[t_k,t_{k+1}]\1n\times\1n\dbR^n,\\
\ns\ds\Th^k(T,\xi,x)=\psi(t_k,\xi,x),\ \ \Th^k(t_j,\xi,x)=\Th^k(t_j+0,\xi,x), \ \ j=N-1,\cdots,k+1,\ \  x\in\dbR^n.\ea\right.\ee
We recall the definition of $\t^\Pi(t)$ (see (\ref{t})), and define
\bel{bar t}\bar\t^\Pi(t)=\sum_{k=0}^{N-2}t_{k+1}I_{[t_k,t_{k+1})}(t)+t_NI_{[t_{N-1},T]}(t),\qq t\in[0,T].\ee
Then
$$0\les\bar\t^\Pi(t)-t\les\|\Pi\|,\qq\forall t\in[0,T],$$
and
$$[\t^\Pi(t),\bar\t^\Pi(t))=[t_k,t_{k+1}),\qq\forall t\in[t_k,t_{k+1}),~0\les k\les N-1.$$
Let
$$\Th^\Pi(t,s,\xi,x)=\sum_{k=0}^{N-1}\Th^k(s,\xi,x)I_{[t_k,t_{k+1})}(t).\qq0\les t\les s\les T,\q x,\xi\in\dbR^n.$$
For $t\in[t_k,t_{k+1})$, with $k=0,1,\cds,N-1$, $s\in[\t^\Pi(t),T]=[t_k,T]$,
$$\ba{ll}
\ns\ds\Th^k(s,\xi,x)=\Th^\Pi(t,s,\xi,x),\q\Th^k_s(s,\xi,x)=\Th^\Pi_s(t,s,\xi,x),\\
\ns\ds\Th^k_x(s,\xi,x)=\Th^\Pi_x(t,s,\xi,x),\q\Th^k_{xx}(s,\xi,x)=\Th^\Pi_{xx}(t,s,\xi,x),\ea$$
and
$$\sum_{\ell=k+1}^{N-1}\Th^\ell(s,x,x)I_{[t_\ell,t_{\ell+1})}(s)=\Th^\Pi(s,s,x,x),\qq s\in[t_{k+1},T].$$
Then the above PDE \eqref{Th(k)} can be written as
\bel{Th(Pi)}\left\{\2n\ba{ll}
\ns\ds\Th^\Pi_s(t,s,\xi,x)+{1\over2}\,\si(s,x)^\top\Th^\Pi_{xx}(t,s,\xi,x)
\si(s,x)+\Th^\Pi_x(t,s,\xi,x)b(s,x)\\
\ns\ds\q+g\big(\t^\Pi(t),s,\xi,x,\Th^\Pi(s,s,x,x),\Th^\Pi_x(t,s,\xi,x)\si(s,x)
\big)=0,\q(s,\xi,x)\in[\bar\t^\Pi(t),T]\times\dbR^n\times\dbR^n,\\
\ns\ds\Th^\Pi_s(t,s,\xi,x)+{1\over2}\,\si(s,x)^\top\Th^\Pi_{xx}(t,s,\xi,x)
\si(s,x)+\Th^\Pi_x(t,s,\xi,x)b(s,x)\\
\ns\ds\q+g\big(\t^\Pi(t),s,\xi,x,\Th^\Pi(s,s,\xi,x),\Th^\Pi_x(t,s,\xi,x)\si(s,x)
\big)=0,\q(s,\xi,x)\in[\t^\Pi(t),\bar\t^\Pi(t))\1n\times\1n\dbR^n\1n\times
\1n\dbR^n,\\
\ns\ds\Th^\Pi(t,T,\xi,x)=\psi(\t^\Pi(t),\xi,x),\qq(t,\xi,x)\in[t_k,T]\times\dbR^n\times\dbR^n,\ea\right.\ee
%
%
%
%
%
%
%
Also,
$$\ba{ll}
\ns\ds\dbE\[\sup_{s\in[t_{k+1},T]}|\bar Y^k(s)-Y^k(s)|^2+\int_{t_{k+1}}^T
|\bar Z^k(s)-Z^k(s)|^2ds\]\\
\ns\ds\les K_1\dbE\(\sum_{\ell=k+1}^{N-1}\int_{t_\ell}^{t_{\ell+1}}\big|g(t_k,s,X(t_k),\Th^\ell(s,X(s),X(s)),
Z^k(s)\big)\\
\ns\ds\qq\qq-g\big(t_k,s,X(t_k),\Th^\ell(s,X(t_\ell),X(s)),Z^k(s)\big)\big|ds\)^2\\
\ns\ds\les K_1L^2\dbE\(\sum_{\ell=k+1}^{N-1}\int_{t_\ell}^{t_{\ell+1}}\big|\Th^\ell(s,X(s),X(s))
-\Th^\ell(s,X(t_\ell),X(s))\big|ds\)^2\\
\ns\ds\les K_0L^2\|\Th_\xi\|_\infty^2\dbE\(\sum_{\ell=k+1}^{N-1}\int_{t_\ell}^{t_{\ell+1}}
\big|X(s)-X(t_\ell)\big|ds\)^2\\
\ns\ds\les K_1L^2\|\Th_\xi\|_\infty^2(T-t_{k+1})\(\sum_{\ell=k+1}^{N-1}\int_{t_\ell}^{t_{\ell+1}}
\dbE|X(s)-X(t_\ell)|^2ds\)\\
\ns\ds\les K_1L^2\|\Th_\xi\|_\infty^2(T-t_{k+1})K_1\(\sum_{\ell=k+1}^{N-1}\int_{t_\ell}^{t_{\ell+1}}
(s-t_\ell)ds\)\\
\ns\ds\les{K_0K_1L^2\|\Th_\xi\|_\infty^2\over2}(T-t_{k+1})\sum_{\ell=k+1}^{N-1}(t_{\ell+1}
-t_\ell)^2\les{K_0K_1L^2\|\Th_\xi\|_\infty^2\over2}(T-t_{k+1})^2\|\Pi\|,\ea$$
where
$$\|\Th_\xi\|_\infty=\max_{k+1\les\ell\les N}\|\Th^\ell_\xi\|_\infty,$$
and
\bel{3.27}\ba{ll}
\ns\ds\dbE\[\sup_{s\in[t_k,t_{k+1})}|\bar Y^k(s)-Y^k(s)|^2
+\int_{t_k}^{t_{k+1}}|\bar Z^k(s)-Z^k(s)|^2ds\]\\
\ns\ds\les K_1\dbE|\bar Y^k(t_{k+1})-Y^k(t_{k+1})|^2
\les{K_0K_1^2L^2\|\Th_\xi\|_\infty^2\over2}(T-t_{k+1})^2\|\Pi\|.\ea\ee
Now, let
$$\left\{\2n\ba{ll}
\ns\ds\bar Y^\Pi(s)=\sum_{k=0}^{N-1}\bar Y^k(s)I_{[t_k,t_{k+1})}(s),\qq s\in[0,T),\\
\ns\ds\bar Z^\Pi(t,s)=\sum_{k=0}^{N-1}\bar Z^k(s)I_{[t_k,t_{k+1})}(t),\qq0\les t\les s\les T.\ea\right.$$
Then
$$\left\{\2n\ba{ll}
\ns\ds\bar Y^\Pi(s)=\sum_{k=0}^{N-1}\Th^k(s,X(t_k),X(s))I_{[t_k,t_{k+1})}(s)
=\Th^\Pi\big(s,s,X(\t^\Pi(s)),X(s)\big),\\
\ns\ds\bar Z^\Pi(t,s)=\sum_{k=0}^{N-1}\Th^k_x(s,X(t_k),X(s))\si(s,X(s))I_{[t_k,t_{k+1})}(t)
=\Th^\Pi\big(t,s,X(\t^\Pi(t)),X(s)\big)\si(s,X(s)).\ea\right.$$
Consequently, for any $s\in[0,T)$, let $s\in[t_k,t_{k+1})$.
$$\ba{ll}
\ns\ds\dbE|Y^\Pi(s)-\Th^\Pi(s,s,X(s),X(s))|^2\les2\dbE|Y^\Pi(s)-\bar Y^\Pi(s)|^2+2\dbE|\bar Y^\Pi(s)-\Th^\Pi(s,s,X(s),X(s))|^2\\ [1mm]
\ns\ds=2\dbE|Y^k(s)-\bar Y^k(s)|^2+2\dbE|\Th^k(s,X(t_k),X(s))-\Th^k(s,X(s),X(s))|^2\\ [1mm]
\ns\ds\les K_0K_1^2L^2\|\Th_\xi\|_\infty^2(T-t_{k+1})^2\|\Pi\|+2\|\Th_\xi\|_\infty^2
\dbE|X(t_k)-X(s)|^2\\ [1mm]
\ns\ds\les K_0K_1^2L^2\|\Th_\xi\|_\infty^2T^2\|\Pi\|
+2\|\Th_\xi\|_\infty^2K_0\|\Pi\|\les K\|\Pi\|.\ea$$
Also,
$$\ba{ll}
\ns\ds\dbE\int_0^T\int_t^T|Z^\Pi(t,s)-\Th^\Pi_x(t,s,X(\t^\Pi(t)),X(s))\si(s,X(s))|^2ds\\
\ns\ds\les\sum_{k=0}^{N-1}\dbE\int_{t_k}^{t_{k+1}}\int_{t_k}^T|Z^k(s)
-\Th^k_x(s,X(s),X(s))\si(s,X(s))|^2dsdt\\
\ns\ds\les\sum_{k=0}^{N-1}\dbE\int_{t_k}^{t_{k+1}}\int_{t_k}^T\(|Z^k(s)-\bar Z^k(s)|^2+|\Th^k_x(s,X(t_k),X(s))-\Th^k_x(s,X(s),X(s))|^2\)dsdt\\
\ns\ds\les\sum_{k=0}^{N-1}2\dbE\int_{t_k}^{t_{k+1}}K\|\Pi\|dt
+\sum_{k=0}^{N-1}2\dbE\int_{t_k}^{t_{k+1}}K\|\Pi\|dt\les K\|\Pi\|.\ea$$
Hence, at the limit (as $\|\Pi\|\to0$), we have the following representation:
\bel{Rep**}\left\{\2n\ba{ll}
\ns\ds Y(s)=\Th(s,s,X(s),X(s)),\qq\qq s\in[0,T],\\
\ns\ds Z(t,s)=\Th_x(t,s,X(t),X(s))\si(s,X(s)),\qq(t,s)\in\D[0,T],\ea\right.\ee
if $\Th(\cd\,,\cd\,,\cd\,,\cd)$ satisfies \rf{PDE0} which is rewritten here:
\bel{PDE**}\left\{\2n\ba{ll}
\ns\ds\Th_s(t,s,\xi,x)+{1\over2}\,\si(s,x)^\top\Th_{xx}(t,s,\xi,x)\si(s,x)
+\Th_x(t,s,\xi,x)b(s,x)\\
\ns\ds\qq+g(t,s,\xi,x,\Th(s,s,x,x),\Th_x(t,s,\xi,x)\si(s,x)\big)=0,\qq(t,s,\xi,x)
\in\D[0,T]\times\dbR^n\times\dbR^n,\\
\ns\ds\Th(t,T,\xi,x)=\psi(t,\xi,x),\qq(t,\xi,x)\in[0,T]\times\dbR^n\times\dbR^n.\ea\right.\ee
The above derivation tells us that if everything is fine, (\ref{Rep**})--(\ref{PDE**}) should give us the right representation. This can actually be proved directly.

\bt{Theorem 3.1.} \sl Let {\rm(H1)--(H2)} hold. Let $\Th:\D[0,T]\times\dbR^n\times\dbR^n\to\dbR$ be the classical solution of the {\it representation} PDE $(\ref{PDE**})$. Let $(Y(\cd),Z(\cd\,,\cd))$ be the adapted solution to the Type-I BSVIE $(\ref{BSVIE1**})$ with $X(\cd)$ being the solution to SDE $(\ref{SDE2})$. Then representation $(\ref{Rep**})$ holds.

\et

\it Proof. \rm For fixed $t\in[0,T)$, applying It\^o's formula to $s\mapsto\Th(t,s,X(t),X(s))$ on $[t,T]$, we have
\bel{}\ba{ll}
\ns\ds d\Th(t,s,X(t),X(s))=\(\Th_s(t,s,X(t),X(s))+\Th_x(t,s,X(t),X(s))b(s,X(s))\\
\ns\ds
\qq\q+{1\over2}\si(s,X(s))^{\top}\Th_{xx}(t,s,X(t),X(s))\si(s,X(s))\)ds
+\Th_x(t,s,X(t),X(s))\si(s,X(s))dW(s).\ea\ee
Since $\Th(\cd\,,\cd\,,\cd\,,\cd)$ satisfies PDE (\ref{PDE**}), one has
\bel{}\ba{ll}
\ns\ds d\Th(t,s,X(t),X(s))=
-g\big(t,s,X(t),X(s),\Th(s,s,X(s),X(s)),\Th_x(t,s,X(t),X(s)\si(s,X(s))\big)ds\\
\ns\ds\qq\qq+\Th_x(t,s,X(t),X(s))\si(s,X(s))dW(s),\ea\ee
and
$$\Th(t,T,X(t),X(T))=\psi(t,X(t),X(T)).$$
Now, we define
\bel{}\ba{ll}
\ns\ds \l(t,s):=\Th(t,s,X(t),X(s)),\ \ Z(t,s):=\Th_x(t,s,X(t),X(s))\si(s,X(s)),\q s\ges t.\ea\ee
Then
\bel{}\ba{ll}
\ns\ds \l(t,s)=\psi(t,X(t),X(T))+\int_s^Tg
\big(t,r,X(t),X(r),\l(r,r),Z(t,r)\big)dr-\int_s^TZ(t,r)dW(r).
\ea\ee
Let $t=s$ and $Y(s):=\l(s,s)$, we then see that $(Y(\cd),Z(\cd\,,\cd))$ satisfies BSVIE (\ref{BSVIE1**}) and desired representation is obtained. \endpf

\section{Representation of Adapted M-solutions for Type-II BSVIEs}

In this section, we are going to establish a representation of adapted M-solutions for Type-II BSVIE \eqref{BSVIE2*}, where both $Z(t,s)$ and $Z(s,t)$ appear in the drift. We still let $X(\cd)$ be the solution to FSDE \eqref{SDE2}. Let us first present the following result which is interesting itself.

\ms

\bp{Proposition 4.1.} \sl Let $\L:[0,T]\times\dbR^n\to\dbR^m$ be continuous. Let {\rm(H1)} hold and the following PDE system admit a classical solution $\G(\cd\,,\cd\,,\cd)$:
\bel{G0}\left\{\2n\ba{ll}
\ns\ds\G_s(t,s,x)+{1\over2}\,\si(s,x)^\top\G_{xx}(t,s,x)\si(s,x)+\G_x(t,s,x)b(s,x)=0,\qq0\les s\les t\les T,\q x\in\dbR^n,\\
\ns\ds\G(t,t,x)=\L(t,x),\qq(t,x)\in[0,T]\times\dbR^n,\ea\right.\ee
where the meaning of $\si(s,x)^\top\G_{xx}(t,s,x)\si(s,x)$ is similar to
{\rm(\ref{si Th si})}. Then
\bel{Rep0}\L(t,X(t))=\dbE\L(t,X(t))+\int_0^t\G_x(t,s,X(s))\si(s,X(s))dW(s),\qq t\in[0,T],\ee

\ep

\it Proof. \rm We consider the following (decoupled) FBSDE on $[0,t]$:
\bel{FBSDE0}\left\{\2n\ba{ll}
\ns\ds dX(s)=b(s,X(s))ds+\si(s,X(s))dW(s),\qq s\in[0,t],\\
\ns\ds d\eta(t,s)=\z(t,s)dW(s),\qq s\in[0,t],\\
\ns\ds\eta(t,t)=\L(t,X(t)),\ea\right.\ee
where $t\in[0,T)$ is a parameter. Then the following representation holds:
\bel{Rep4.4}\left\{\2n\ba{ll}
\ns\ds\eta(t,s)=\G(t,s,X(s)),\\
\ns\ds\z(t,s)=\G_x(t,s,X(s))\si(s,X(s)),\ea\right.\qq s\in[0,t],\ee
where $\G(t,\cd\,,\cd)$ is the solution to \rf{G0}. Consequently,
$$\L(t,X(t))=\eta(t,t)=\eta(t,0)+\int_0^t\z(t,s)dW(s).$$
Taking expectation, we have
$$\dbE\L(t,X(t))=\eta(t,0).$$
Therefore, \rf{Rep0} follows. \endpf

\ms

From the above, we see that when $(t,s)\mapsto\G_x(t,s,x)$ and $s\mapsto\si(s,x)$ are continuous, the map $(t,s)\mapsto\z(t,s)$ is continuous (see \rf{Rep4.4}).

\ms

Now, we consider Type-II BSVIE \eqref{BSVIE2*}. Let $(Y(\cd),Z(\cd\,,\cd))$ be the adapted M-solution. Then we have \eqref{M}. Suppose
$$Y(t)=\L(t,X(t)),\qq t\in[0,T],$$
for some undetermined continuous function $\L(\cd\,,\cd)$. By Proposition \ref{Proposition 4.1.}, we have
$$Z(t,s)=\G_x(t,s,X(s))\si(s,X(s)),\qq0\les s\les t\les T.$$
Thus, switching $s$ and $t$, one has
$$Z(s,t)=\G_x(s,t,X(t))\si(t,X(t)),\qq0\les t\les s\les T.$$
We consider the following Type-I BSVIE:
\bel{BSVIE4}\ba{ll}
\ns\ds Y(t)=\psi(t,X(t),X(T))+\int_t^Tg\big(t,s,X(t),X(s),Y(s),Z(t,s),\G_x(s,t,X(t))\si(t,X(t))
\big)ds\\
\ns\ds\qq\qq\qq\qq\qq\qq\qq\qq\qq\qq-\int_t^TZ(t,s)dW(s),\qq t\in[0,T].\ea\ee
If we let
$$\wt g(t,s,\xi,x,y,z)=g\big(t,s,\xi,x,y,z,\G_x(s,t,\xi)\si(t,\xi)\big),$$
then \eqref{BSVIE4} becomes
\bel{BSVIE5}\ba{ll}
\ns\ds Y(t)=\psi(t,X(t),X(T))+\int_t^T\wt g\big(t,s,X(t),X(s),Y(s),Z(t,s)\big)ds-\int_t^TZ(t,s)dW(s),\q t\in[0,T].\ea\ee
Now, from the result of the previous section, we have the following representation:
\bel{Rep**-M}\left\{\2n\ba{ll}
\ns\ds Y(s)=\Th(s,s,X(s),X(s)),\qq s\in[0,T],\\ [2mm]
\ns\ds Z(t,s)=\Th_x(t,s,X(t),X(s))\si(s,X(s)),\qq0\les t\les s\les T,\\ [2mm]
\ns\ds Z(t,s)=\G_x(t,s,X(s))\si(s,X(s)),\qq0\les s\les t\les T,\ea\right.\ee
with $(\G,\Th)$ being the solution to \rf{PDE**-M} which is rewritten here
\bel{PDE**-M*}\left\{\2n\ba{ll}
\ns\ds\G_s(t,s,x)+{1\over2}\,\si(s,x)^\top\G_{xx}(t,s,x)\si(s,x)+\G_x(t,s,x)b(s,x)=0,\q0\les s\les t\les T,~x\in\dbR^n,\\
\ns\ds\Th_s(t,s,\xi,x)+{1\over2}\,\si(s,x)^\top\Th_{xx}(t,s,\xi,x)\si(s,x)+\Th_x(t,s,\xi,x)b(s,x)\\
\ns\ds\qq+g(t,s,\xi,x,\Th(s,s,x,x),\Th_x(t,s,\xi,x)\si(s,x),\G_x(s,t,\xi)\si(s,x)\big)=0,\\
\ns\ds\qq\qq\qq\qq\qq\qq\qq\qq\qq\qq\qq(t,s,\xi,x)\in\D[0,T]\times\dbR^n\times\dbR^n,\\
\ns\ds\G(t,t,x)=\Th(t,t,x,x),\qq(t,x)\in[0,T]\times\dbR^n,\\
\ns\ds\Th(t,T,\xi,x)=\psi(t,\xi,x),\qq(t,\xi,x)\in[0,T]\times\dbR^n\times\dbR^n.\ea\right.\ee

\ms

We now state the representation result as follows.

\bt{Theorem 4.2.} \sl Let $\Th,\ \G$ be the solution of representation PDE $(\ref{PDE**-M*})$. Let $(Y(\cd),Z(\cd\,,\cd))$ be the adapted M-solution to Type-II BSVIE $(\ref{BSVIE2*})$ with $X(\cd)$ being the solution to SDE $(\ref{SDE2})$. Then representation $(\ref{Rep**-M})$ holds.

\et

\it Proof. \rm For given $t\in[0,T)$, applying It\^o's formula to $s\mapsto\Th(t,s,X(t),X(s))$ on $[t,T]$, one has
\bel{}\ba{ll}
\ns\ds d\Th(t,s,X(t),X(s))=\(\Th_s(t,s,X(t),X(s))+\Th_x(t,s,X(t),X(s))b(s,X(s))\\
\ns\ds
\qq\q+{1\over2}\si(s,X(s))^{\top}\Th_{xx}(t,s,X(t),X(s))\si(s,X(s))\)ds
+\Th_x(t,s,X(t),X(s))\si(s,X(s))dW(s).\ea\ee
Since $\Th$ satisfies the second PDE of (\ref{PDE**-M*}), one has
\bel{}\ba{ll}
\ns\ds d\Th(t,s,X(t),X(s))=\Th_x(t,s,X(t),X(s))\si(s,X(s))dW(s)\\
\ns\ds
-g\big(t,s,X(t),X(s),\Th(s,s,X(s),X(s)),\Th_x(t,s,X(t),X(s)\si(s,X(s)),
\Gamma(s,t,X(t))\si(s,X(s))\big)ds,\ea\ee
and
$$\Th(t,T,X(t),X(T))=\psi(t,X(t),X(T)).$$
Set
\bel{}\ba{ll}
\ns\ds \l(t,s):=\Th(t,s,X(t),X(s)),\ \ Z(t,s):=\Th_x(t,s,X(t),X(s))\si(s,X(s)),\q s\ges t.
\ea\ee
Then
\bel{}\ba{ll}
\ns\ds \l(t,s)=\psi(t,X(t),X(T))+\int_s^Tg
\big(t,r,X(t),X(r),\l(r,r),Z(t,r),\Gamma(r,t,X(t))\si(r,X(r))\big)dr\\
\ns\ds\qq\qq\qq\qq\qq\qq-\int_s^TZ(t,r)dW(r).
\ea\ee
Let $t=s$ and $Y(s):=\l(s,s)$, we then obtain
\bel{}\ba{ll}
\ns\ds Y(t)=\psi(t,X(t),X(T))+\int_t^Tg
\big(t,r,X(t),X(r),\l(r,r),Z(t,r),\Gamma(r,t,X(t))\si(r,X(r))\big)dr\\
\ns\ds\qq\qq\qq\qq\qq-\int_t^TZ(t,r)dW(r).\ea\ee
Note that
$$Y(t)=\Th(t,t,X(t),X(t))=\G(t,t,X(t))$$
where $\G$ satisfies the first PDE in (\ref{PDE**-M*}). By Proposition \ref{Proposition 4.1.}, we know that
$$Y(t)=\dbE Y(t)+\int_0^t\G_x(t,s,X(s))\si(s,X(s))dW(s),\q t\ges s.$$
Consequently, by defining $Z(t,s):=\G_x(t,s,X(s))\si(s,X(s))$ with $t\ges s,$ we can rewrite above BSVIE as
\bel{}\ba{ll}
\ns\ds Y(t)=\psi(t,X(t),X(T))+\int_t^Tg
\big(t,r,X(t),X(r),\l(r,r),Z(t,r),Z(r,t)\big)dr-\int_t^TZ(t,r)dW(r).
\ea\ee
The conclusion then follows easily. \endpf

\section{Well-posedness of the Representation PDEs}

In this section, we will establish the well-posedness of the representation PDEs \rf{PDE**} (which is a copy of \rf{PDE0}) and \rf{PDE**-M*} (which is a copy of \rf{PDE**-M}), in certain sense. Let us first look at the representation PDE \rf{PDE**} for Type-I BSVIEs, which is recalled here, for convenience:
\bel{PDE5.1}\left\{\2n\ba{ll}
\ds\Th_s(t,s,\xi,x)+{1\over2}\si(s,x)^\top\Th_{xx}(t,s,\xi,x)\si(s,x)+\Th_x(t,s,\xi,x)b(s,x)\\
\ns\ds\qq\qq\qq+g(t,s,\xi,x,\Th(s,s,x,x),\Th_x(t,s,\xi,x)\si(s,x)\big)=0,
\q(t,s,\xi,x)\1n\in\1n\D[0,T]\1n\times\1n\dbR^n\1n\times\1n\dbR^n,\\
\ns\ds\Th(t,T,\xi,x)=\psi(t,\xi,x),\qq(t,\xi,x)\in[0,T]\times\dbR^n\times\dbR^n.\ea\right.\ee
If we denote $\Th=(\Th^1,\cds,\Th^m)$ and
$${1\over2}\si(s,x)\si(s,x)^\top=a(s,x)=\big(a_{ij}(s,x)\big),\q b(s,x)=(b_1(s,x),\cds,b_n(s,x))^\top,$$
then \rf{PDE5.1} can be rewritten as the following system (parameterized by $(t,\xi)\in[0,T)\times\dbR^n$):
\bel{PDE-system}\left\{\2n\ba{ll}
\ds\Th^k_s(t,s,\xi,x)+\sum_{i,j=1}^na_{ij}(s,x)\Th^k_{x_ix_j}(t,s,\xi,x)
+\sum_{i=1}^nb_i(s,x)\Th^k_{x_i}(t,s,\xi,x)\\
\ns\ds\qq\qq+g^k(t,s,\xi,x,\Th(s,s,x,x),\Th_x(t,s,\xi,x)\si(s,x)\big)=0,\qq(s,x)
\in[t,T]\times\dbR^n,\\
\ns\ds\Th^k(t,T,\xi,x)=\psi^k(t,\xi,x),\qq x\in\dbR^n,\qq1\les k\les m,\ea\right.\ee
which is a quasilinear parabolic system for unknown functions $\Th^1,\cds,\Th^m$, with the same leading part for each equation.

\subsection{Linear parabolic PDEs}

To study parabolic system \rf{PDE5.1} or its equivalent form \rf{PDE-system}, let us first adopt some notations from \cite{Ladyzenskaja 1968} (Chapter 1, pp.7--8). For any suitable function $\f:[S,T]\times\dbR^n\to\dbR$, with $\a\in(0,1)$ and $S\in[0,T)$, let
\bel{|f|}\left\{\2n\ba{ll}
\ns\ds|\f|^{(0)}=\|\f\|_{L^\infty([S,T]\times\dbR^n)},\q|\f|^{(1)}=|\f|^{(0)}
+|\f_x|^{(0)}_,\q|\f|^{(2)}=|\f|^{(1)}+|\f_s|^{(0)}+|\f_{xx}|^{(0)},\\ [2mm]
\ns\ds\lan\f\ran{}\1n_s^{({\a\over2})}=\3n\sup_{{s,s'\in[S,T],x\in\dbR^n}
\atop{s\ne s'}}\3n{|\f(s,x)-\f(s',x)|\over|s-s'|^{\a\over2}},\qq
\lan\f\ran{}\1n^{(\a)}_x=\3n\sup_{{s\in[S,T],x,x'\in\dbR^n}
\atop{0<|x-x'|\les1}}\3n{|\f(s,x)-\f(s,x')|\over|x-x'|^\a},\\ [2mm]
\ns\ds\lan\f\ran\1n{}^{(\a)}=\lan\f\ran{}\1n_s^{({\a\over2})}
+\lan\f\ran\1n{}_x^{(\a)},\qq|\f|^{(\a)}=|\f|^{(0)}+\lan\f\ran{}\1n^{(\a)},\\ [1mm]
\ns\ds|\f|^{(1+\a)}=|\f|^{(1)}+\lan\f_x\ran{}\1n^{(\a)}+\lan\f\ran{}\1n_s^{({1+\a\over2})},\\ [1mm]
\ns\ds|\f|^{(2+\a)}=|\f|^{(2)}+\lan\f_s\ran{}\1n^{(\a)}+\lan\f_{xx}\ran{}\1n^{(\a)}
+\lan\f_x\ran\1n{}^{({1+\a\over2})}_s.\ea\right.\ee
When $[S,T]\times\dbR^n$ needs to be emphasized, we use, say, $|\f|^{(1)}_{[S,T]\times\dbR^n}$, etc. We denote
$$C^{{\a\over2},\a}([S,T]\times\dbR^n)=\Big\{\f:[S,T]\times\dbR^n\to\dbR\bigm||\f|^{(\a)}_{[S,T]
\times\dbR^n}<\infty\Big\}.$$
Clearly, $\f(\cd\,,\cd)\in C^{{\a\over2},\a}([S,T]\times\dbR^n)$ if and only if $\f(\cd\,,\cd)\in L^\infty([S,T]\times\dbR^n)$ and
$$|\f(s,x)-\f(s',x')|\les K\(|s-s'|^{\a\over2}+|x-x'|^\a\),\qq\forall s,s'\in[S,T],~x,x'\in\dbR^n,~|x-x'|\les1.$$
Also, we denote
$$\ba{ll}
\ns\ds C^{{1+\a\over2},1+\a}([S,T]\times\dbR^n)=\Big\{\f:[S,T]\to\dbR^n\bigm|
|\f|^{(1+\a)}_{[S,T]\times\dbR^n}<\infty\Big\},\\
\ns\ds C^{1+{\a\over2},2+\a}([S,T]\times\dbR^n)=\Big\{\f:[S,T]\to\dbR^n\bigm|
|\f|^{(2+\a)}_{[S,T]\times\dbR^n}<\infty\Big\}.\ea$$

\ms

Let us consider the following Cauchy problem for linear equation:
\bel{linear parabolic PDE}\left\{\2n\ba{ll}
\ds v_s(t,s,\xi,x)+\sum_{i,j=1}^na_{ij}(s,x)v_{x_ix_j}(s,x)
+\sum_{i=1}^nb_i(s,x)v_{x_i}(s,x)+f(s,x)=0,\q(s,x)
\in[t,T]\times\dbR^n,\\
\ns\ds v(T,x)=h(x),\qq x\in\dbR^n.\ea\right.\ee
We introduce the following hypotheses:

\ms

{\bf(P1)} Operator $\cL$ is {\it uniformly parabolic}, i.e., there exist constants $\bar\l_0>0$ such that
$$\sum_{i,j=1}^na_{ij}(s,x)\xi_i\xi_j\equiv\lan a(s,x)\xi,\xi\ran\ges\bar\l_0|\xi|^2,\qq\forall\xi\in\dbR^n,\q(s,x)\in[0,T]\times\dbR^n.$$

{\bf(P2)} Functions $a_{ij}(s,x),b_i(s,x)$ are continuous and bounded, and for some $\a\in(0,1)$,
$$\left\{\2n\ba{ll}
\ds|a_{ij}(s,x)-a_{ij}(s',x')|\les K\big(|s-s'|^{\a\over2}+|x-x'|^\a\big),\\
\ns\ds|b_i(s,x)-b_i(s,x')|\les K|x-x'|^\a,\ea\right.\qq\q(s,x),(s',x')\in[0,T]\times\dbR^n.$$

\ms

By \cite{Ladyzenskaja 1968} (Chapter IV, Sections 13--14) (see also \cite{Friedman 1964}, Chapter 1, Section 7) we have the following result.

\bp{Friedman} \sl Let {\rm(P1)--(P2)} hold. Assume that
$$f(\cd\,,\cd)\in C^{{\a\over2},\a}([0,T]\times\dbR^n),\q h(\cd)\in C^{2+\a}(\dbR^n),$$
for some $\a\in(0,1)$. Then Cauchy problem \rf{linear parabolic PDE} admits a unique classical solution $v(\cd\,,\cd)\in C^{1+{\a\over2},2+\a}([0,T]\times\dbR^n)$. Moreover, $v(\cd\,,\cd)$ is represented as follows:
$$v(s,x)=\int_{\dbR^n}G(s,x;T,\eta)h(\eta)d\eta+\int_s^T\2n\int_{\dbR^n}G(s,x;\t,\eta)f(\t,\eta)d\eta d\t,\qq(s,x)\in[0,T]\times\dbR^n.$$
Here $G(s,x;\t,\eta)$ is called the {\it fundamental solution} of the parabolic operator $\cL$, which satisfies the following: There exists a $\l>0$ such that for any $x,\eta\in\dbR^n$,
\bel{|G|<}\left\{\2n\ba{ll}
\ds|G(s,x;\t,\eta)|\les{K\over(\t-s)^{n\over2}}e^{-\l{|\eta-x|^2\over\t-s}},\\
\ns\ds|G_x(s,x;\t,\eta)|\les{K\over(\t-s)^{n+1\over2}}e^{-\l|\eta-x|^2\over\t-s},\\
\ns\ds|G_s(s,x;\t,\eta)|+|G_{xx}(s,x;\t,\eta)|\les{K\over(\t-s)^{n+2\over2}}e^{-\l|\eta-x|^2\over\t-s},\ea\right.\qq0\les s<\t\les T.\ee
%
%
%
%
%
%
%
%
Moreover,
\bel{|v|(2+a)}|v|^{(2+\a)}_{[0,T]\times\dbR^n}\les K\big(|f|^{(\a)}_{[0,T]\times\dbR^n}+|h|^{(2+\a)}_{\dbR^n}\big),\ee
and for any $S\in[0,T)$,
\bel{|v|(1+a)}|v|^{(1+\a)}_{[S,T]\times\dbR^n}\les|h|^{(1+\a)}_{ \dbR^n}+ K(T-S)^{\a\over2}\big(|f|^{(\a)}_{[S,T]\times\dbR^n}
+|h|^{(2+\a)}_{\dbR^n}\big).\ee

\ep

\it Proof. \rm The proof up to \rf{|v|(2+a)} is standard (see \cite{Ladyzenskaja 1968,Friedman 1964}). Let us
look at \rf{|v|(1+a)} which will play a crucial role below.

\ms

Note that by defining
$$\ba{ll}
\ns\ds\wt v(s,x)=v(s,x)-h(x),\\
\ns\ds\wt f(s,x)=f(s,x)+{1\over2}\si(s,x)^\top h_{xx}(x)+h_x(x)b(s,x),\qq(s,x)\in[0,T]\times\dbR^n,\ea$$
we see that $\wt v(\cd\,,\cd)$ is the solution to the following:
\bel{parabolic PDE*}\left\{\2n\ba{ll}
\ds\wt v_s(s,x)+{1\over2}\si(s,x)^\top\wt v_{xx}(s,x)\si(s,x)+\wt v_x(s,x)b(s,x)+\wt f(s,x)=0,\qq(s,x)\in[0,T]\times\dbR^n,\\
\ns\ds\wt v(T,x)=0,\qq x\in\dbR^n.\ea\right.\ee
Thus, the following representation holds:
\bel{}\wt v(s,x)=\int_s^T\3n\int_{\dbR^n}G(s,x;\t,\eta)\wt f(\t,\eta)d\eta d\t,\qq(s,x)\in[0,T]\times\dbR^n.\ee
%
%
%
%
%
Following the steps of proving the inequality \rf{|v|(2+a)} in \cite{Ladyzenskaja 1968} (Chapter IV, Section 14), we have the following useful estimates:
\bel{5.9}\left\{\2n\ba{ll}
\ds|\wt v|^{(0)}_{[S,T]\times\dbR^n}\les K(T-S)|\wt f|^{(\a)}_{[S,T]\times\dbR^n},\\ [2mm]
\ns\ds|\wt v_x|^{(0)}_{[S,T]\times\dbR^n}\les K(T-S)^{1+\a\over2}|\wt f|^{(\a)}_{[S,T]\times\dbR^n},\\ [2mm]
\ns\ds|\wt v_s|^{(0)}_{[S,T]\times\dbR^n}+|\wt v_{xx}|^{(0)}_{[S,T]\times\dbR^n}\les K(T-S)^{\a\over2}|\wt f|^{(\a)}_{[S,T]\times\dbR^n}.\ea\right.\ee
Also, for any $s,s'\in[S,T]$, $x,x'\in\dbR^n$,
\bel{5.10}\left\{\2n\ba{ll}
\ds|\wt v_x(s,x)-\wt v_x(s',x)|\les K|s-s'|^{1+\a\over2}|\wt f|^{(\a)}_{[S,T]\times\dbR^n},\\
\ns\ds|\wt v_{xx}(s,x)-\wt v_{xx}(s',x)|\les K|s-s'|^{\a\over2}|\wt f|^{(\a)}_{[S,T]\times\dbR^n},\ea\right.\qq(s,x),(s',x)\in[S,T]\times\dbR^n.\ee
Now, from \eqref{5.9}, for any $s,s'\in[S,T]$, $x,x'\in\dbR^n$,
$|x-x'|\les1$, we further have
$$\ba{ll}
\ns\ds|\wt v(s,x)-\wt v(s',x)|\les\int_0^1|\wt v_s(s'+\m(s-s'),x)|d\m\,|s-s'|\les
K(T-S)^{\a\over2}|\wt f|^{(\a)}_{[S,T]\times\dbR^n}|s-s'|\\
\ns\ds\qq\qq\qq\qq\les K|\wt f|^{(\a)}_{[S,T]\times\dbR^n}
\[(T-S)|s-s'|^{\a\over2}\]\land\[(T-S)^{1\over2}|s-s'|^{1+\a\over2}\],\\
\ns\ds|\wt v(s,x)-\wt v(s,x')|\les\int_0^1|\wt v_x(s,x'+\m(x-x'))|d\m\,|x-x'|\les K(T-S)^{1+\a\over2}|\wt f|^{(\a)}_{[S,T]\times\dbR^n}|x-x'|\\
\ns\ds\qq\qq\qq\qq\les K(T-S)^{1+\a\over2}|\wt f|^{(\a)}_{[S,T]\times\dbR^n}|x-x'|^\a,\ea$$
which leads to
\bel{5.4}\lan\wt v\ran\1n{}^{(\a)}_{[S,T]\times\dbR^n}
\equiv\lan\wt v\ran\1n{}^{({\a\over2})}_{s,[S,T]\times\dbR^n}
+\lan\wt v\ran\1n{}^{(\a)}_{x,[S,T]\times\dbR^n}\les
K(T-S)^{1+\a\over2}|\wt f|^{(\a)}_{[S,T]\times\dbR^n}.\ee
Next, the first inequality in \rf{5.10} implies that
$$|\wt v_x(s,x)-\wt v_x(s',x)|\les K|\wt f|^{(\a)}_{[S,T]\times\dbR^n}
|s-s'|^{1+\a\over2}\les K(T-S)^{1\over2}|\wt f|^{(\a)}_{[S,T]\times\dbR^n}
|s-s'|^{\a\over2}.$$
%
%
Similar to the above, for any $s\in[S,T]$, $x,x'\in\dbR^n$, $|x-x'|\les1$, making use of the third inequality in \rf{5.9}, we have
$$\ba{ll}
\ns\ds|\wt v_x(s,x)-\wt v_x(s,x')|\les\int_0^1|\bar v_{xx}(s,x'+\m(x-x'))d\m|\,|x-x'|\les
K(T-S)^{\a\over2}|\wt f|^{(\a)}_{[S,T]\times\dbR^n}|x-x'|,\ea$$
which leads to
\bel{5.5}\lan\wt v_x\ran\1n{}^{(\a)}_{[S,T]\times\dbR^n}\les K(T-S)^{\a\over2}|\wt f|^{(\a)}_{[S,T]\times\dbR^n}.\ee
Hence, combining the above, we end up with
\bel{5.6}\ba{ll}
\ns\ds|\wt v|^{(1+\a)}_{[S,T]\times\dbR^n}\1n\equiv|\wt v|^{(0)}_{[S,T]\times\dbR^n}
+|\wt v_x|^{(0)}_{[S,T]\times\dbR^n}+\lan\wt v\ran\1n{}^{({1+\a\over2})}_{s,[S,T]\times\dbR^n}+\lan\wt v_x\ran\1n{}^{(\a)}_{[S,T]\times\dbR^n}\les K(T-S)^{\a\over2}|\wt f|^{(\a)}_{[S,T]\times\dbR^n}.\ea\ee
This implies \rf{|v|(1+a)}. \endpf

\ms

\subsection{The first representation PDE}

Now, let us return to the first representation PDE \eqref{PDE5.1}. We impose the following further assumption.

\ms

{\bf(H3)} The maps $b(s,x)$, $\si(s,x)$, $\psi(t,\xi,x)$, and $g(t,s,\xi,x,y,z)$ are bounded, have all required differentiability with bounded derivatives. Moreover, there exists a constant $\bar\si>0$ such that
\bel{si>0}|\si(s,x)\xi|\ges\bar\si|\xi|^2,\qq\forall(s,x,\xi)\in[0,T]\times\dbR^n\times\dbR^n.\ee

\ms

The above assumption is much more than enough. However, in this paper, we prefer not to get into the most generality in this aspect, to reduce the complexity of presentation. Also, we will extend $(t,s)\mapsto g(t,s,\xi,x,y,z)$ from $\D[0,T]$ to $[0,T]$ by letting
$$g(t,s,\xi,x,y,z)=g(s,t,\xi,x,y,z),\qq\forall t,s\in[0,T].$$
We now state the following theorem.

\bt{1st} \sl Let {\rm(H3)} hold. Then for any $(t,\xi)\in[0,T)\times\dbR^n$, system \rf{PDE5.1} admits a unique classical solution $\Th(t,\cd\,,\xi,\cd)\in C^{1+{\a\over2},2+\a}([t,T]\times\dbR^n)$, and the following holds
\bel{|Th|}\sup_{(t,\xi)\in[0,T]\times\dbR^n}\big|\Th(t,\cd\,,\xi,\cd)\big|^{(2+\a)}_{[t,T]\times\dbR^n}\les K
\(1+\sup_{(t,\xi)\in[0,T]\times\dbR^n}|\psi(t,\xi,\cd)|^{(2+\a)}_{\dbR^n}\).\ee
Moreover, if $\h\Th(t,\cd\,,\xi,\cd)$ is the solution to the system \eqref{PDE5.1} with the pair $(\psi,g)$ replaced by $(\h\psi,\h g)$ that also satisfies {\rm(H3)}, then the following stability estimate holds:
\bel{|Th-Th|}\ba{ll}
\ns\ds\sup_{(t,\xi)\in[0,T]\times\dbR^n}\big|\Th(t,\cd\,,\xi,\cd)-\h\Th(t,\cd\,,\xi,\cd)\big|^{(2+\a)}_{[t,T]\times\dbR^n}\\
\ns\ds\les K\sup_{(t,\xi)\in[0,T]\times\dbR^n}
 \(|\psi(t,\xi,\cd)-\h\psi(t,\xi,\cd)|^{(2+\a)}_{\dbR^n}+|g(t,\cd\,,\xi,\cd)-\h g(t,\cd\,,\xi,\cd)|^{(\a)}_{[t,T]\times\dbR^n)}\),\\
\ns\ds\qq\qq\qq\qq\qq\qq\qq\qq\qq\qq\forall(t,\xi)\in[0,T)\times\dbR^n,\ea\ee
%
%
%
where
\bel{gg}\ba{ll}
\ns\ds g(t,s,\xi,x)=g\big(t,s,\xi,x,\Th(s,s,x,x),\Th_x(t,s,\xi,x)
\si(s,x)\big),\\
\ns\ds\h g(t,s,\xi,x)=\h g\big(t,x,\xi,x,\Th(s,s,x,x),\Th_x(t,s,\xi,x)
\si(s,x)\big),\ea\q(t,s,\xi,x)\in\D[0,T]\times\dbR^n\times\dbR^n.\ee

\et

\it Proof. \rm We split the proof into several steps.

\ms

\it Step 1. A reduction. \rm First, we let
$$\wt\Th(t,s,\xi,x)=\Th(t,s,\xi,x)-\psi(t,\xi,x),\qq(t,s)\in\D[0,T],~
x,\xi\in\dbR^n.$$
Then $\Th(\cd\,,\cd\,,\cd\,,\cd)$ is a solution of \rf{PDE5.1} if
and only if $\wt\Th(\cd\,,\cd\,,\cd\,,\cd)$ is a solution to the following:
\bel{PDE-system-wt-Th}\left\{\2n\ba{ll}
\ds\wt\Th_s(t,s,\xi,x)+{1\over2}\si(s,x)^\top\wt\Th_{xx}(t,s,\xi,x)
\si(s,x)+\wt\Th_x(t,s,\xi,x)b(s,x)\\ [1mm]
\ns\ds\qq+\wt g(t,s,\xi,x,\wt\Th(s,s,x,x),\wt\Th_x(t,s,\xi,x)
\si(s,x)\big)=0,\qq(s,x)\in[t,T]\times\dbR^n,\\ [1mm]
\ns\ds\wt\Th(t,T,\xi,x)=0,\qq x\in\dbR^n,\ea\right.\ee
where
$$\ba{ll}
\ns\ds\wt g(t,s,\xi,x,y,z)=g(t,s,\xi,x,y+\psi(t,\xi,x),
\psi_x(t,\xi,x)\si(s,x)+z\si(s,x)\big)\\
\ns\ds\qq\qq\qq\qq+{1\over2}\si(s,x)^\top\psi_{xx}(t,\xi,x)\si(s,x)
+\psi_x(t,\xi,x)b(s,x)=0,\q(s,x)\in[t,T]\times\dbR^n,\ea$$
Hence, without loss of generality, we may consider \rf{PDE5.1} with $\psi(t,\xi,x)\equiv0$.

\ms

\it Step 2. The solution map of a parabolic PDE. \rm Let
\bel{hD}\h\D[S,T]=\([0,S]\times[S,T]\)\bigcup\D[S,T]\equiv\Big\{(t,s)\in[0,T]^2\bigm|
0\les t\vee S\les s\les T\Big\}.\ee
See the figure below for illustration.


\setlength{\unitlength}{.01in}
~~~~~~~~~~~~~~~~~~~~~~~~~~~~~~~~~~~~~~~~~~\begin{picture}(290,270)
\put(0,0){\vector(1,0){240}}
\put(0,0){\vector(0,1){240}}
\put(150,0){\line(0,1){150}}
\put(200,0){\line(0,1){200}}
\put(0,50){\line(1,0){50}}
\put(0,100){\line(1,0){100}}
\put(0,150){\color{blue}\line(1,0){150}}
\put(0,200){\line(1,0){200}}
\thicklines
\put(0,0){\color{red}\line(1,1){200}}
\put(-25,45){\makebox(0,0)[b]{$T-3\d$}}
\put(65,125){\makebox(0,0){$\h\D[S',S]$}}
\put(70,175){\makebox(0,0){$\h\D[S,T]$}}
\put(-40,95){\makebox(0,0)[b]{$S'=T-2\d$}}
\put(-35,145){\makebox(0,0)[b]{$S=T-\d$}}
\put(-10,195){\makebox(0,0)[b]{$T$}}
\put(150,-15){\makebox(0,0)[b]{$S$}}
\put(200,-15){\makebox(0,0)[b]{$T$}}
\put(250,-5){\makebox{$t$}}
\put(0,245){\makebox{$s$}}
\end{picture}

\bs

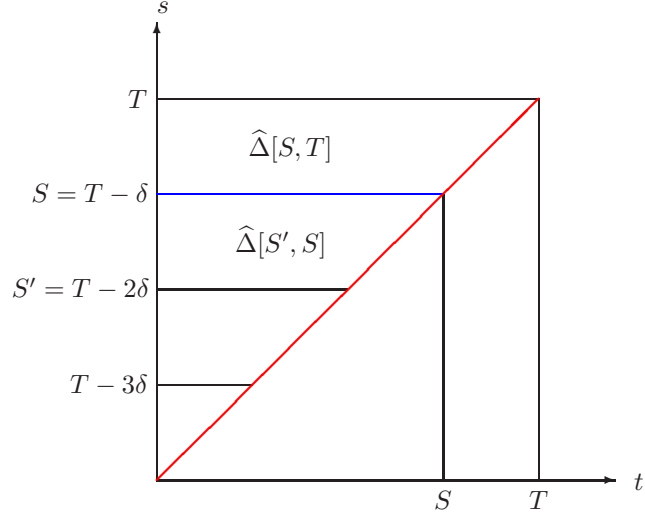
\centerline{(Figure 2)}

\bs

Let $\sX[S,T]$ be the set of all measurable functions $\th:\h\D[S,T]\times\dbR^{2n}\to\dbR^m$ such that
\bel{|th|}\ba{ll}
\ns\ds\|\th\|_{\sX[S,T]}\equiv\sup_{(t,\xi)\in[0,T]\times\dbR^n}
\(|\th(t,\cd\,,\xi,\cd)|^{(1+\a)}_{[t\vee S,T]\times\dbR^n}
+|\th_t(t,\cd\,,\xi,\cd)|^{(0)}_{[t\vee S,T]\times\dbR^n}
+|\th_\xi(t,\cd\,,\xi,\cd)|^{(0)}_{[t\vee S,T]\times\dbR^n}\\
\ns\ds\qq\qq\qq\qq\qq\qq\qq\qq+|\th_{xt}(t,\cd\,,\xi,\cd)|^{(0)}_{[t\vee S,T]\times\dbR^n}
+|\th_{x\xi}(t,\cd\,,\xi,\cd)|^{(0)}_{[t\vee S,T]\times\dbR^n}\)\\
\ns\ds\equiv\sup_{(t,\xi)\in[0,T]\times\dbR^n}\(|\th(t,\cd\,,\xi,\cd)|^{(0)}_{[t\vee S,T]\times\dbR^n}
+|\th_\eta(t,\cd\,,\xi,\cd)|^{(0)}_{[t\vee S,T]\times\dbR^n}
+\lan\th_\eta(t,\cd\,,\xi,\cd)\ran\1n{}^{(\a)}_{[t\vee S,T]\times\dbR^n}\\
\ns\ds\qq\qq\qq\qq+\lan\th(t,\cd\,,\xi,\cd)\ran\1n{}^{({1+\a\over2})}_{s,
[t\vee S,T]\times\dbR^n}+|\th_t(t,\cd\,,\xi,\cd)|^{(0)}_{[t\vee S,T]\times\dbR^n}+|\th_\xi(t,\cd\,,\xi,\cd)|^{(0)}_{[t\vee S,T]\times\dbR^n}\\ [2mm]
\ns\ds\qq\qq\qq\qq+|\th_{xt}(t,\cd\,,\xi,\cd)|^{(0)}_{[t\vee S,T]\times\dbR^n}
+|\th_{x\xi}(t,\cd\,,\xi,\cd)|^{(0)}_{[t\vee S,T]\times\dbR^n}\)<\infty.\ea\ee
Clearly, $\|\cd\|_{\sX[S,T]}$ is a norm under which $\sX[S,T]$ is a Banach space.

\ms

Let $S\in[0,T)$ be fixed. For any $\th(\cd\,,\cd\,,\cd\,,\cd)
\in\sX[S,T]$, denote
\bel{g}g(t,\t,\xi,\eta)=g\big(t,\t,\xi,\eta,\th(\t,\t,\eta,\eta),\th_\eta(t,\t,\xi,\eta)
\si(\t,\eta)\big).\ee
We claim that
\bel{|g|a}|g(t,\cd\,,\xi,\cd)|^{(\a)}_{[t\vee S,T]\times\dbR^n}\les K\big(1+\|\th\|_{\sX[S,T]}\big),\qq\forall(t,\xi)\in[S,T]\times\dbR^n.\ee
Let us prove this. By boundedness of $g$, one has
\bel{|g|<K}\big|g\big(t,\t,\xi,\eta,\th(\t,\t,\eta,\eta),\th_\eta(t,\t,\xi,\eta)
\si(\t,\eta)\big)\big|\les K.\ee
Next, for $\t,\t'\in[t\vee S,T]$ and $\xi,\eta\in\dbR^n$,
\bel{|g|t}\ba{ll}
\ns\ds\big|g(t,\t,\xi,\eta,\th(\t,\t,\eta,\eta),\th_\eta(t,\t,\xi,\eta)\si(\t,\eta))
-g(t,\t',\xi,\eta,\th(\t',\t',\eta,\eta),\th_\eta(t,\t',\xi,\eta)\si(\t',\eta))
\big|\\ [2mm]
\ns\ds\les K\(|\t-\t'|+|\th(\t,\t,\eta,\eta)-\th(\t',\t',\eta,\eta)|+|\th_\eta(t,\t,\xi,\eta)\si(\t,\eta)
-\th_\eta(t,\t',\xi,\eta)\si(\t',\eta)|\)\\ [3mm]
\ns\ds\les K\(|\t-\t'|+\sup_{(\bar t,\bar\xi)\in[t\vee S,T]\times\dbR^n}|\th_t(\bar t,\cd\,,\bar\xi,\cd)|^{(0)}_{[t\vee S,T]\times\dbR^n}|\t-\t'|\\
\ns\ds\qq\qq+\sup_{(\bar t,\bar\xi)\in[t\vee S,T]\times\dbR^n}\lan\th(\bar t,\cd\,,\bar\xi,\cd)\ran\1n
{}^{({\a\over2})}_{\t,[t\vee S,T]\times\dbR^n}|\t-\t'|^{\a\over2}\\ [4mm]
\ns\ds\qq\qq+|\th_\eta(t,\cd\,,\xi,\cd)|^{(0)}_{[t\vee S,T]\times\dbR^n}|\t-\t'|+\lan\th_\eta(t,
\cd\,,\xi,\cd)\ran\1n{}^{({\a\over2})}_{\t,[t\vee S,T]\times\dbR^n}|\t-\t'|^{\a\over2}\)\\
\ns\ds\les K\(1+\|\th\|_{\sX[S,T]}\)|\t-\t'|^{\a\over2}.\ea\ee
Likewise, for $\t\in[t\vee S,T]$ and $\xi,\eta,\eta'\in\dbR^n$,
\bel{|g|eta}\ba{ll}
\ns\ds|g(t,\t,\xi,\eta,\th(\t,\t,\eta,\eta),\th_\eta(t,\t,\xi,\eta)\si(\t,\eta))
-g(t,\t,\xi,\eta',\th(\t,\t,\eta',\eta'),\th_\eta(t,\t,\xi,\eta')\si(\t,\eta'))|\\
\ns\ds\les K\(|\eta-\eta'|+|\th(\t,\t,\eta,\eta)
-\th(\t,\t,\eta',\eta')|+|\th_\eta(t,\t,\xi,\eta)\si(\t,\eta)
-\th_\eta(t,\t,\xi,\eta')\si(\t,\eta')|\)\\
\ns\ds\les K\(|\eta-\eta'|+\sup_{(\bar t,\bar\xi)\in[t\vee S,T]\times\dbR^n}|\th_\xi(\bar t,\cd\,,\bar\xi,\cd)|^{(0)}_{[t\vee S,T]\times\dbR^n}|\eta-\eta'|\\
\ns\ds\qq\qq+\sup_{(\bar t,\bar\xi)\in[t\vee S,T]\times\dbR^n}\lan\th(\bar t,\cd\,,\bar\xi,\cd)\ran\1n{}^{(\a)}_{[t\vee S,T]\times\dbR^n}|\eta-\eta'|^\a\\
\ns\ds\qq\qq+|\th_\eta(t,\cd\,,\xi,\cd)|^{(0)}_{[t\vee S,T]\times\dbR^n}|\eta-\eta'|
+\lan\th_\eta(t,\cd\,,\xi,\cd)\ran\1n{}^{(\a)}_{\eta,[t\vee S,T]\times\dbR^n}|\eta-\eta'|^\a\)\\
\ns\ds\les K\(1+\|\th\|_{\sX[S,T]}\)|\eta-\eta'|^\a.\ea\ee
Combining \rf{|g|<K}--\rf{|g|eta}, we obtain \rf{|g|a}. Now, for any $\th(\cd\,,\cd\,,\cd\,,\cd)\in\sX[S,T]$, consider the following linear parabolic system, parametrized by $(t,\xi)\in[0,T)\times\dbR^n$:
\bel{PDE5.35}\left\{\2n\ba{ll}
\ns\ds\Th_s(t,s,\xi,x)+{1\over2}\si(s,x)^\top\Th_{xx}(t,s,\xi,x)
\si(s,x)+\Th_x(t,s,\xi,x)b(s,x)\\ [1mm]
\ns\ds\qq\qq\qq+g(t,s,\xi,x,\th(s,s,x,x),\th_x(t,s,\xi,x)\si(s,x)\big)=0,
\qq(s,x)\in[t\vee S,T]\times\dbR^n,\\ [1mm]
\ns\ds\Th(t,T,\xi,x)=0,\qq x\in\dbR^n.\ea\right.\ee
Then the corresponding solution $\Th(t,\cd\,,\xi,\cd)$ uniquely
exists and the following holds:
\bel{Th=int}\Th(t,s,\xi,x)=\int_s^T\3n\int_{\dbR^n}
G(s,x;\t,\eta)g\big(t,\t,\xi,\eta,\th(\t,\t,\eta,\eta),
\th_\eta(t,\t,\xi,\eta)\si(\t,\eta)\big)d\eta d\t.\ee
Due to \rf{|g|a}, we have $\Th(t,\cd\,,\xi,\cd)\in C^{2+\a}([t\vee S,T]\times\dbR^n)$. On the other hand, by \rf{|v|(1+a)} and \rf{|g|a}, we have
\bel{5.25}\ba{ll}
\ns\ds|\Th(t,\cd\,,\xi,\cd)|^{(1+\a)}_{[t\vee S,T]\times\dbR^n}\les
K(T-S)^{\a\over2}|g(t,\cd\,,\xi,\cd)|^{(\a)}_{[t\vee S,T]\times\dbR^n}\les K(T-S)^{\a\over2}(1+\|\th\|_{\sX[S,T]}).\ea\ee
Next,
$$\ba{ll}
\ns\ds\Th_t(t,s,\xi,x)=\int_s^T\3n\int_{\dbR^n}
G(s,x;\t,\eta)g_t\big(t,\t,\xi,\eta,\th(\t,\t,\eta,\eta),
\th_\eta(t,\t,\xi,\eta)\si(\t,\eta)\big)d\eta d\t\\
\ns\ds\qq\qq\qq+\int_s^T\3n\int_{\dbR^n}
G(s,x;\t,\eta)g_z\big(t,\t,\xi,\eta,\th(\t,\t,\eta,\eta),
\th_\eta(t,\t,\xi,\eta)\si(\t,\eta)\big)\th_{\eta t}(t,\t,\xi,\eta)
\si(\t,\eta)d\eta d\t,\ea$$
and
$$\ba{ll}
\ns\ds\Th_{xt}(t,s,\xi,x)=\int_s^T\3n\int_{\dbR^n}
G_x(s,x;\t,\eta)g_t\big(t,\t,\xi,\eta,\th(\t,\t,\eta,\eta),
\th_\eta(t,\t,\xi,\eta)\si(\t,\eta)\big)d\eta d\t\\
\ns\ds\qq\qq\qq+\int_s^T\3n\int_{\dbR^n}
G_x(s,x;\t,\eta)g_z\big(t,\t,\xi,\eta,\th(\t,\t,\eta,\eta),
\th_\eta(t,\t,\xi,\eta)\si(\t,\eta)\big)\th_{\eta t}(t,\t,\xi,\eta)
\si(\t,\eta)d\eta d\t.\ea$$
Using \rf{|G|<}, one has
$$\ba{ll}
\ns\ds|\Th_{xt}(t,s,\xi,x)|\les K\int_s^T\3n\int_{\dbR^n}
|G_x(s,x;\t,\eta)|d\eta d\t+K\int_s^T\3n\int_{\dbR^n}
|G_x(s,x;\t,\eta)|\,|\th_{\eta t}(t,\t;\xi,\eta)|d\eta d\t\\
\ns\ds\les K\int_s^T\3n\int_{\dbR^n}{1\over(\t-s)^{n+1\over2}}
e^{-\l|\eta-x|^2\over\t-s}d\eta d\t+K|\th_{xt}|^{(0)}_{\h\D[S,T]
\times\dbR^{2n}}\int_s^T\3n\int_{\dbR^n}{1\over(\t-s)^{n+1\over2}}
e^{-\l|\eta-x|^2\over\t-s}d\eta d\t\\
\ns\ds\les K(T-S)^{1\over2}+K(T-s)^{1\over2}
|\th_{xt}|^{(0)}_{\h\D[S,T]\times\dbR^{2n}},\ea$$
where
$$|\th_{xt}|^{(0)}_{\h\D[S,T]\times\dbR^{2n}}=\sup_{(t,s)\in\h\D[S,T],\xi,x
\in\dbR^n}|\th_{xt}(t,s,\xi,x)|.$$
Similarly,
$$\ba{ll}
\ns\ds|\Th_t(t,s,\xi,x)|\les K\int_s^T\3n\int_{\dbR^n}
|G(s,x;\t,\eta)|d\eta d\t+K\int_s^T\3n\int_{\dbR^n}
|G(s,x;\t,\eta)|\,|\th_{\eta t}(t,\t;\xi,\eta)|d\eta d\t\\ [2mm]
\ns\ds\les K\int_s^T\3n\int_{\dbR^n}{1\over(\t-s)^{n\over2}}
e^{-\l|\eta-x|^2\over\t-s}d\eta d\t+K|\th_{xt}|^{(0)}_{\h\D[S,T]
\times\dbR^{2n}}\int_s^T\3n\int_{\dbR^n}{1\over(\t-s)^{n\over2}}
e^{-\l|\eta-x|^2\over\t-s}d\eta d\t\\ [4mm]
\ns\ds\les K(T-s)+K(T-s)|\th_{xt}|^{(0)}_{\h\D[S,T]\times\dbR^{2n}},\ea$$
Consequently, it holds that
\bel{|Th_t|}|\Th_t|^{(0)}_{\h\D[S,T]\times\dbR^{2n}}+|\Th_{xt}|^{(0)}_{\h
\D[S,T]\times\dbR^{2n}}\les K(T-S)^{1\over2}\(1+|\th_{xt}|^{(0)}_{\h\D[S,T]\times\dbR^{2n}}\).\ee
Likewise, we have
\bel{|Th_xi|}|\Th_\xi|^{(0)}_{\h\D[S,T]\times\dbR^{2n}}
+|\Th_{x\xi}|^{(0)}_{\h\D[S,T]\times\dbR^{2n}}\les K(T-S)^{1\over2}\(1+|\th_{x\xi}|^{(0)}_{\h\D[S,T]\times\dbR^{2n}}\).\ee
Then combining \rf{5.25} with \rf{|Th_t|}--\rf{|Th_xi|}, one obtains
\bel{|Th|<}\|\Th\|_{\sX[S,T]}\les K(T-S)^{\a\over2}\big(1+\|\th\|_{\sX[S,T]}\big),\ee
with $K>0$ being an absolute constant. Hence, we have defined a map $\cS:\sX[S,T]\to\sX[S,T]$ by
$$\cS[\th(\cd\,,\cd\,,\cd\,,\cd)]=\Th(\cd\,,\cd\,,\cd\,,\cd),\qq\forall\th(\cd\,,\cd\,,\cd\,,\cd)\in\sX[S,T].$$
Moreover, we shrink $T-S>0$ (if necessary) so that $K(T-S)^{\a\over2}\les{1\over2}$. Then for any $M\ges1$,
$$\|\Th\|_{\sX[S,T]}\les{1\over2}\big(1+\|\th\|_{\sX[S,T]}\big)\les M,
\qq\forall\|\th\|_{\sX[S,T]}\les M.$$
Thus, $\cS$ maps a ball in $\sX[S,T]$, centered at 0 with radius $M$ to itself.

\ms

\it Step 3. Contraction of the solution map. \rm Let $\th,\h\th\in\sX[S,T]$ such that
$$\|\th\|_{\sX[S,T]},\;\|\h\th\|_{\sX[S,T]}\les M,$$
with $S\in[0,T)$ and $M>0$ obtained as above. Let
$$\ba{ll}
\ns\ds\Th(t,s,\xi,x)=\cS[\th](t,s,\xi,x),\qq\h\Th(t,s,\xi,x)=\cS[\h\th](t,s,\xi,x),\\
\ns\ds v(t,s,\xi,x)=\Th(t,s,\xi,x)-\h\Th(t,s,\xi,x).\ea$$
Then $v(t,\cd\,,\xi,\cd)$ satisfies the following:
\bel{PDE5.9}\left\{\2n\ba{ll}
\ns\ds v_s(t,s,\xi,x)+{1\over2}\si(s,x)^\top v_{xx}(t,s,\xi,x)\si(s,x)+v_x(t,s,\xi,x)b(s,x)+f(t,s,\xi,x)=0,\\
\ns\ds\qq\qq\qq\qq\qq\qq\qq\qq\qq\qq(s,x)\in[t\vee S,T]\times\dbR^n,\\
\ns\ds v(t,T,\xi,x)=0,\qq x\in\dbR^n,\ea\right.\ee
with
$$f(t,s,\xi,x)=g(t,s,\xi,x,\th(s,s,x,x),\th_x(t,s,\xi,x)\si(s,x))
-g(t,s,\xi,x,\h\th(s,s,x,x),\h\th_x(t,s,\xi,x)\si(s,x)).$$
Consequently,
\bel{5.33}v(t,s,\xi,x)=\int_s^T\2n\int_{\dbR^n}G(s,x;\t,\eta)f(t,\t,\xi,\eta)d\eta d\t.\ee
Similar to the proof of \rf{|g|a}, we can show that
%
\bel{|f|a}|f(t,\cd\,,\xi,\cd)|^{(\a)}_{[t\vee S,T]\times\dbR^n}\les K\|\th-\h\th\|_{\sX[S,T]},\qq\forall(t,\xi)\in[S,T]\times\dbR^n.\ee
Then one has
\bel{|Th-Th|}\|\Th-\h\Th\|_{\sX[S,T]}\les  K(T-S)^{\a\over2}\|\th-\h\th\|_{\sX[S,T]},\ee
for some absolute constant $K>0$. Hence, by choosing $\d=T-S>0$ small, we obtain a contraction mapping $\cS$ on $\sX[T-\d,T]$. Consequently, $\cS$ has a unique fixed point which is the solution $\Th(\cd\,,\cd\,,\cd\,,\cd)$ on $\h\D[T-\d,T]\times\dbR^n\times\dbR^n$, which is rewritten here:
\bel{PDE5.40}\left\{\2n\ba{ll}
\ns\ds\Th_s(t,s,\xi,x)+{1\over2}\si(s,x)^\top\Th_{xx}(t,s,\xi,x)
\si(s,x)+\Th_x(t,s,\xi,x)b(s,x)\\ [1mm]
\ns\ds\qq\qq\qq+g(t,s,\xi,x,\Th(s,s,x,x),\Th_x(t,s,\xi,x)\si(s,x)\big)=0,
\qq(t,s)\in\h\D[S,T]\times\dbR^n\times\dbR^n,\\ [1mm]
\ns\ds\Th(t,T,\xi,x)=0,\qq(t,\xi,x)\in[0,T]\times\dbR^n\times\dbR^n.\ea\right.\ee
Note that $(t,\xi)\mapsto g(t,s,\xi,x,y,z)$ is assumed to be continuously differentiable with bounded derivatives. Therefore, by a standard argument, we know that
$$|\Th_t(t,\cd\,,\xi,\cd)|^{(2+\a)}_{[t\vee (T-\d),T]\times\dbR^n}+|\Th_\xi(t,\cd\,,\xi,\cd)|^{(2+\a)}_{[t\vee (T-\d),T]\times\dbR^n}<\infty.$$
Now, we denote
$$\bar\psi(t,\xi,x)=\Th(t,T-\d,\xi,x),\qq(t,\xi,x)\in[0,T-\d]\times
\dbR^n\times\dbR^n,$$
and consider the following equation:
\bel{PDE5.40}\left\{\2n\ba{ll}
\ns\ds\Th_s(t,s,\xi,x)+{1\over2}\si(s,x)^\top\Th_{xx}(t,s,\xi,x)
\si(s,x)+\Th_x(t,s,\xi,x)b(s,x)\\ [1mm]
\ns\ds\qq\qq\qq+g(t,s,\xi,x,\Th(s,s,x,x),\Th_x(t,s,\xi,x)\si(s,x)\big)=0,
\qq(s,x)\in[t,T-\d]\times\dbR^n,\\ [1mm]
\ns\ds\Th(t,T-\d,\xi,x)=\bar\psi(t,\xi,x),\qq x\in\dbR^n.\ea\right.\ee
Then, we may repeat the above procedure, to get a unique solution $\Th(t,\cd\,,\xi,\cd)$ on $\h\D[T-2\d,T-\d]\times\dbR^n\times\dbR^n$.
By continuing such a procedure, we obtain the existence and unique solution $\Th(t,\cd\,,\xi,\cd)$ to the representation PDE \eqref{PDE5.1}, and \eqref{|Th|} holds.

\ms

\it Step 4. Stability estimates. \rm Let $(\h\psi,\h g)$ be another pair of maps such that, together with $b$ and $\si$, satisfy (H3) as well. Let $\h\Th$ be the corresponding solution. Let
$$v(t,s,\xi,x)=\Th(t,s,\xi,x)-\h\Th(t,s,\xi,x).$$
Then $v(t,\cd\,,\xi,\cd)$ satisfies the following:
\bel{PDE5.9*}\left\{\2n\ba{ll}
\ns\ds v_s(t,s,\xi,x)+{1\over2}\si(s,x)^\top v_{xx}(t,s,\xi,x)\si(s,x)+v_x(t,s,\xi,x)b(s,x)+f(t,s,\xi,x)+\h f(t,s,\xi,x)=0,\\
\ns\ds\qq\qq\qq\qq\qq\qq\qq\qq\qq\qq\qq\qq\qq(s,x)\in[t,T]\times\dbR^n,\\
\ns\ds v(t,T,\xi,x)=\f(t,\xi,x),\qq x\in\dbR^n,\ea\right.\ee
with
$$\ba{ll}
\ns\ds f(t,s,\xi,x)=\h g(t,s,\xi,x,\Th(s,s,x,x),\Th_x(t,s,\xi,x)\si(s,x))
-\h g(t,s,\xi,x,\h\Th(s,s,x,x),\h\Th_x(t,s,\xi,x)\si(s,x)),\\
\ns\ds\bar f(t,s,\xi,x)=g(t,s,\xi,x,\Th(s,s,x,x),\Th_x(t,s,\xi,x)\si(s,x))
-\h g(t,s,\xi,x,\Th(s,s,x,x),\Th_x(t,s,\xi,x)\si(s,x)),\\
\ns\ds\f(t,\xi,x)=\psi(t,\xi,x)-\h\psi(t,\xi,x).\ea$$
Then
\bel{5.39}\ba{ll}
\ns\ds v(t,s,\xi,x)=\int_{\dbR^n}G(s,x;T,\eta)\f(t,\xi,\eta)d\eta+\int_s^T\int_{\dbR^n}G(s,x;\t,\eta)
\bar f(t,\t,\xi,\eta)d\eta d\t\\ [2mm]
\ns\ds\qq\qq\qq\qq+\int_s^T\int_{\dbR^n}G(s,x;\t,\eta)f(t,\t,\xi,\eta)d\eta d\t.\ea\ee
From (\ref{|v|(1+a)}), one has
$$\ba{ll}
\ns\ds\Big|\int_{\dbR^n}G(\cd\,,\cd\,;T,\eta)\f(t,\xi,\eta)d\eta+\int_t^T\2n\int_{\dbR^n}G(\cd\,,\cd\,;
\t,\eta)\bar f(t,\t,\xi,\eta)d\eta d\t\Big|^{(2+\a)}_{[t,T]\times\dbR^n}\\
\ns\ds\les K\(|\f(t,\xi,\cd)|^{(2+\a)}_{\dbR^n}+|\bar f(t,\cd\,,\xi,\cd)|^{(\a)}_{[t,T]\times\dbR^n}\).\ea$$
Also, by Step 3 above, we have
$$\Big\|\int_\cd^T\int_{\dbR^n}G(\cd\,,\cd;\t,\eta)f(\cd\,,\t,\cd\,,\eta)
d\eta d\t\Big\|_{\sX[S,T]}
\les K(1+M)(T-S)^{\a\over2}\|v\|_{\sX[S,T]}.$$
Hence, for $T-S>0$ small, we obtain
$$\|v\|_{\sX[S,T]}\les K\sup_{(t,\xi)\in[S,T]\times\dbR^n}\(|\f(t,\xi,\cd)|^{(2+\a)}_{\dbR^n}+|\bar f(t,\cd\,,\xi,\cd)|^{(\a)}_{[t,T]\times\dbR^n}\).$$
Repeating the same argument, we obtain
$$\|v\|_{\sX[0,T]}\les K\sup_{(t,\xi)\in[0,T]\times\dbR^n}\(|\f(t,\xi,\cd)|^{(2+\a)}_{\dbR^n}+|\bar f(t,\cd\,,\xi,\cd)|^{(\a)}_{[0,T]\times\dbR^n}\).$$
Then \eqref{|Th-Th|} follows. \endpf

\ms

\subsection{The second representation PDE}

We now look at the second representation PDE. Similar to the previous subsection, without loss of generality, we again assume that  $\psi(t,\xi,x)=0$. Thus, we consider the following family of parabolic systems (parameterized by $(t,\xi)\in(S,T)\times\dbR^n$):
\bel{PDE-[0,T]}\left\{\2n\ba{ll}
\ns\ds\G_s(t,s,x)+{1\over2}\si(s,x)^\top\G_{xx}(t,s,x)\si(s,x)+\G_x(t,s,x)b(s,x)
=0,\qq(s,x)\in[0,t]\times\dbR^n,\\
\ns\ds\Th_s(t,s,\xi,x)+{1\over2}\si(s,x)^\top\Th_{xx}(t,s,\xi,x)\si(s,x)
+\Th_x(t,s,\xi,x)b(s,x)\\
\ns\ds\qq+g(t,s,\xi,x,\Th(s,s,x,x),\Th_x(t,s,\xi,x)\si(s,x),\G_\xi(s,t,\xi)
\si(s,x)\big)=0,\qq(s,x)\in[t,T]\times\dbR^n,\\
\ns\ds\G(t,t,x)=\Th(t,t,x,x),\qq x\in\dbR^n,\\
\ns\ds\Th(t,T,\xi,x)=0,\qq x\in\dbR^n.\ea\right.\ee
Note that for any given $(t,\xi)\in(0,T)\times\dbR^n$, the equation for $\G(t,\cd\,,\cd)$ is to be solved on $[0,t]$ and the equation for $\Th(t,\cd\,,\xi,\cd)$ is to be solved on $[t,T]$. The coupling appears at two places: $\G_\xi(s,t,\xi)$ (with $0\les t\les s$) appears in the equation for $\Th(t,\cd\,,\xi,\cd)$ and $\Th(t,t,x,x)$ appears as the terminal value for $\G(t,s,x)$ at $s=t$.


\setlength{\unitlength}{.01in}
~~~~~~~~~~~~~~~~~~~~~~~~~~~~~~~~~~~~~~~~~~\begin{picture}(290,270)
\put(0,0){\vector(1,0){240}}
\put(0,0){\vector(0,1){240}}
\put(200,0){\line(0,1){200}}
\put(0,200){\line(1,0){200}}
\thicklines
\put(0,0){\color{red}\line(1,1){200}}
\put(65,125){\makebox(0,0){$\G_\xi(s,t,\xi)$}}
\put(65,155){\makebox(0,0){$\Th(t,s,\xi,x)$}}
\put(125,65){\makebox(0,0){$\G(t,s,\xi)$}}
\put(-10,195){\makebox(0,0)[b]{$T$}}
\put(200,-15){\makebox(0,0)[b]{$T$}}
\put(250,-5){\makebox{$t$}}
\put(0,245){\makebox{$s$}}
\end{picture}

\ms

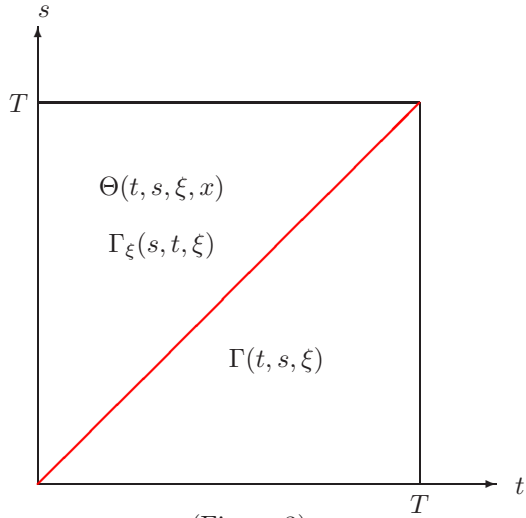
\centerline{(Figure 3)}

\bs

See the above figure for the domains in which $\Th(t,x,\xi,x)$ and $\G(t,s,\xi)$ are defined. Let us make some observations. Suppose $(\G(\cd\,,\cd\,,\cd),\Th(\cd\,,\cd\,,\cd\,,\cd))$ is a solution to the above \rf{PDE-[0,T]}. Then
\bel{G=int}\G(t,s,x)=\int_{\dbR^n}G(s,x;t,\eta)\Th(t,t,\eta,\eta)d\eta,
\qq0\les s\les t\les T,\q x\in\dbR^n,\ee
and thus,
$$\G_\xi(\t,t,\xi)=\int_{\dbR^n}G_\xi(t,\xi;\t,\bar\eta)
\Th(\t,\t,\bar\eta,\bar\eta)d\bar\eta,\qq0\les t\les\t\les T,\q\xi\in\dbR^n.$$
On the other hand,
\bel{Th=int*}\ba{ll}
\ns\ds\Th(t,s,\xi,x)=\int_s^T\3n\int_{\dbR^n}G(s,x;\t,\eta)g\big(t,\t,\xi,\eta,
\Th(\t,\t,\eta,\eta),\Th_\eta(t,\t,\xi,\eta)\si(\t,\eta),\G_\xi(\t,t,\xi)
\si(\t,\eta)\big)d\eta d\t\\ [2mm]
\ns\ds=\int_s^T\3n\int_{\dbR^n}G(s,x;\t,\eta)g\(t,\t,\xi,\eta,
\Th(\t,\t,\eta,\eta),\Th_\eta(t,\t,\xi,\eta)\si(\t,\eta),\\
\ns\ds\qq\qq\qq\qq\[\int_{\dbR^n}
G_\xi(t,\xi;\t,\bar\eta)\Th(\t,\t,\bar\eta,\bar\eta)d\bar\eta
\]\si(\t,\eta)\)d\eta d\t,\q(t,s,\xi,x)\1n\in\1n\D[0,T]\1n\times\1n\dbR^n\1n\times\1n\dbR^n.\ea\ee
The above tells us that if \rf{PDE-[0,T]} admits a classical solution $(\G,\Th)$, then $\Th$ must be a solution to the above nonlinear integral equation \rf{Th=int*}. Conversely, if nonlinear integral equation \rf{Th=int*} admits a smooth solution $\Th$, by defining $\G$ as \rf{G=int}, we have a solution $(\G,\Th)$ to the second representation PDE \rf{PDE-[0,T]}. Hence, we could introduce the following definition.

\bde{mild solution} \rm A pair of functions $(\G,\Th)$ is called a {\it mild solution} of \rf{PDE-[0,T]} if $\Th$ is a solution to the integral equation \rf{Th=int} and $\G$ is defined by \rf{G=int}.

\ede

Now, we introduce the following hypothesis for system \eqref{PDE-[0,T]}.

\ms

{\bf(H4)} The maps $b(s,x)$, $\si(s,x)$, and $g(t,s,\xi,x,y,z,\z)$ are bounded and have all needed orders of bounded derivatives. Moreover \eqref{si>0} holds for some $\bar\si>0$.

\ms

The main result of this subsection is the following theorem.

\bt{mild} \sl Let {\rm (H4)} hold. Then \eqref{PDE-[0,T]} admits a unique mild solution $(\G,\Th)$.

\et

\it Proof. \rm Let $1<p<2$, $S\in[0,T)$ and recall $\h\D[S,T]$ defined by \rf{hD}. Let $\sY[S,T]$ be the set of all functions $\th:\h\D[S,T]\times\dbR^n\times\dbR^n\to\dbR^m$ such that
$$\|\th\|_{\sY[S,T]}=\sup_{t\in[0,T]}\(\int_{t\vee S}^T\sup_{\xi,x\in
\dbR^n}|\th_x(t,s,\xi,x)|^pds\)^{1\over p}+\sup_{{(t,s)\in\h\D[S,T]}\atop{\xi,x\in\dbR^n}}|\th(t,s,\xi,x)|<\infty.$$
Clearly, $\|\cd\|_{\sY[S,T]}$ is a norm under which $\sY[S,T]$ is a Banach space.

\ms

For any $\th\in\sY[S,T]$, define
\bel{}\ba{ll}
\ns\ds\cS[\th](t,s,\xi,x)=\int_s^T\3n\int_{\dbR^n}G(s,x;\t,\eta)g\(t,\t,\xi,\eta,\th(\t,\t,\eta,\eta),
\th_\eta(t,\t,\xi,\eta)\si(\t,\eta),\\
\ns\ds\qq\qq\qq\qq\qq\qq\qq\qq\qq\qq\[\int_{\dbR^n}G_\xi(t,\xi;\t,\bar\eta)\th(\t,\t,\bar\eta,\bar\eta)
d\bar\eta\]\si(\t,\eta)\)d\eta d\t,\\ [2mm]
\ns\ds\qq\qq\qq\qq\qq\qq\qq\qq\qq\qq\qq\qq(t,s,\xi,x)\in\h\D[S,T]\times\dbR^n\times\dbR^n.\ea\ee
Note that
$$\cS[0](t,s,\xi,x)=\int_s^T\3n\int_{\dbR^n}G(s,x;\t,\eta)
g(t,\t,\xi,\eta,0,0,0)d\eta d\t,\q(t,s,\xi,x)\in\h\D[S,T]
\times\dbR^n\times\dbR^n.$$
Thus,
$$\cS[0](s,s,x,x)=\int_s^T\3n\int_{\dbR^n}G(s,x;\t,\eta)
g(s,\t,x,\eta,0,0,0)d\eta d\t,\q(s,x)\in[S,T]\times\dbR^n,$$
and
$$\cS[0]_x(t,s,\xi,x)=\int_s^T\3n\int_{\dbR^n}G_x(s,x;\t,\eta)
g(t,\t,\xi,\eta,0,0,0)d\eta d\t,\q(t,s,\xi,x)\in\h\D[S,T]
\times\dbR^n\times\dbR^n.$$
Then
$$|\cS[0](t,s,\xi,x)|\les K\int_s^T\3n\int_{\dbR^n}{1\over(\t-s)^{n\over2}}e^{-\l{|\eta-x|^2
\over\t-s}}d\eta d\t=K\int_s^T\3n\int_{\dbR^n}e^{-\l z^2}dzd\t\les K.$$
%
%
%
Also,
$$|\cS[0]_x(t,s,\xi,x)|\les K\int_s^T\3n\int_{\dbR^n}{1\over(\t-s)^{n+1\over2}}e^{-\l{|\eta-x|^2
\over\t-s}}d\eta d\t=K\int_s^T\3n\int_{\dbR^n}{1\over(\t-s)^{1\over2}}e^{-\l z^2}dzd\t\les K.$$
Hence, $\cS[0]\in\sY[S,T]$.

\ms

Next, let $\th,\h\th\in\sY[S,T]$. We estimate the following:
$$\ba{ll}
\ns\ds|\cS[\,\th\,](t,s,\xi,x)-\cS[\,\h\th\,](t,s,\xi,x)|\les K\int_s^T\3n\int_{\dbR^n}
|G(s,x;\t,\eta)|\(|\th(\t,\t,\eta,\eta)-\h\th(\t,\t,\eta,\eta)|\\
\ns\ds\qq\qq+|\th_\eta(t,\t,\xi,\eta)-\th_\eta(t,\t,\xi,\eta)|
+\int_{\dbR^n}|G_\xi
(t,\xi,\t,\bar\eta)|\,|\th(\t,\t,\bar\eta,\bar\eta)-\h\th(\t,\t,
\bar\eta,\bar\eta)|d\bar\eta\)d\eta d\t\ea$$
$$\ba{ll}
\ns\ds=K\int_s^T\3n\int_{\dbR^n}|G(s,x;\t,\eta)|\,
|\th_\eta(t,\t,\xi,\eta)-\h\th_\eta(t,\t,\xi,\eta)|d\eta d\t\\
\ns\ds\qq+K\int_s^T\3n\int_{\dbR^n}\(|G(s,x;\t,\eta)|
+|G_\xi(t,\xi;\t,\eta)|\int_{\dbR^n}|G(s,x;\t,\bar\eta)d\bar\eta\)
|\th(\t,\t,\eta,\eta)-\h\th(\t,\t,\eta,\eta)|d\eta d\t\\
\ns\ds\les K\int_s^T\3n\[\int_{\dbR^n}{e^{-\l{|\eta-x|^2\over\t-s}}\over
(\t-s)^{n\over2}}d\eta\]\sup_{x,\eta\in\dbR^n}|\th_\eta(t,\t,\xi,
\eta)-\h\th_\eta(t,\t,\xi,\eta)|d\t\\
\ns\ds\q+K\[\int_s^T\3n\int_{\dbR^n}\({e^{-\l{|\eta-x|^2\over\t-s}}
\over(\t-s)^{n\over2}}+{e^{-\l{|\eta-x|^2\over\t-t}}\over
(\t-t)^{n+1\over2}}\int_{\dbR^n}{e^{-\l{|\bar\eta-x|^2\over\t-s}}
\over(\t-s)^{n\over2}}d\bar\eta\)d\eta d\t\]\sup_{(\t,\eta)\in[s,T]\times\dbR^n}
|\th(\t,\t,\eta,\eta)-\h\th(\t,\t,\eta,\eta)|\\
\ns\ds\les K\int_s^T\[\int_{\dbR^n}e^{-\l z^2}dz\]\sup_{\xi,x\in\dbR^n}|\th_x(s,\t,\xi,x)-\h\th_x(s,\t,\xi,x)|d\t\\
\ns\ds\q+K\[\int_s^T\int_{\dbR^n}e^{-\l z^2}dz+\int_{\dbR^n}{e^{-\l z^2}\over(\t-t)^{1\over2}}\int_{\dbR^n}e^{-\l\bar z^2}d\bar z\)dzd\t\]\sup_{(\t,\eta)\in[s,T]\times\dbR^n}
|\th(\t,\t,\eta,\eta)-\h\th(\t,\t,\eta,\eta)|\\
\ns\ds\les K\int_s^T\sup_{\xi,x\in\dbR^n}|\th_x(s,\t,\xi,x)
-\h\th_x(s,\t,\xi,x)|d\t\\
\ns\ds\qq\qq+K\[(T-t)^{1\over2}-(s-t)^{1\over2}\]\sup_{(\t,\eta)\in[s,T]\times\dbR^n}
|\th(\t,\t,\eta,\eta)-\h\th(\t,\t,\eta,\eta)|\\
\ns\ds\les K(T-s)^{p-1\over p}\(\int_s^T\sup_{\xi,x\in\dbR^n}|\th_x(s,\t,\xi,x)
-\h\th_x(s,\t,\xi,x)|^pd\t\)^{1\over p}\\
\ns\ds\qq+K(T-s)^{1\over2}\sup_{{(t,\t)\in\D[s,T]}\atop{\xi,x\in\dbR^n}}
|\th(t,\t,\xi,x)-\h\th(t,\t,\xi,x)|.\ea$$
This leads to (note ${p-1\over p}<{1\over2}$)
\bel{}\sup_{{(t,s)\in\h\D[S,T]}\atop{\xi,x\in\dbR^n}}\big|\cS[\,\th\,]
(t,s,\xi,x)-\cS[\,\h\th\,](t,s,\xi,x)|\les K(T-S)^{p-1\over p}\|\th-\h\th\|_{\cY[S,T]}.\ee
Next,
$$\ba{ll}
\ns\ds|\cS[\,\th\,]_x(t,s,\xi,x)-\cS[\,\h\th\,]_x(t,s,\xi,x)|\les K\int_s^T\3n\int_{\dbR^n}
|G_x(s,x;\t,\eta)|\(|\th(\t,\t,\eta,\eta)-\h\th(\t,\t,\eta,\eta)|\\
\ns\ds\qq+|\th_\eta(t,\t,\xi,\eta)-\h\th_\eta(t,\t,\xi,\eta)|+\int_{\dbR^n}|G_\xi
(t,\xi,\t,\bar\eta)|\,|\th(\t,\t,\bar\eta,\bar\eta)-\h\th(\t,\t,\bar\eta,
\bar\eta)|d\bar\eta\)
d\eta d\t\\
\ns\ds=K\int_s^T\3n\int_{\dbR^n}|G_x(s,x;\t,\eta)|\,
|\th_\eta(t,\t,\xi,\eta)-\h\th_\eta(t,\t,\xi,\eta)|d\eta d\t\\
\ns\ds\qq+K\int_s^T\3n\int_{\dbR^n}\(|G_x(s,x;\t,\eta)|
+|G_\xi(t,\xi;\t,\eta)|\int_{\dbR^n}|G_x(s,x;\t,\bar\eta)|d\bar\eta\)
|\th(\t,\t,\eta,\eta)-\h\th(\t,\t,\eta,\eta)|d\eta d\t\\
\ns\ds\les K\int_s^T\3n\[\int_{\dbR^n}{e^{-\l{|\eta-x|^2\over\t-s}}\over
(\t-s)^{n+1\over2}}d\eta\]\sup_{\xi,\eta\in\dbR^n}|\th_\eta(t,\t,
\xi,\eta)-\h\th_\eta(t,\t,\xi,\eta)|d\t\\
\ns\ds\q+K\[\int_s^T\3n\int_{\dbR^n}\2n\({e^{-\l{|\eta-x|^2\over\t-s}}
\over(\t-s)^{n+\over2}}\1n+\1n{e^{-\l{|\eta-x|^2\over\t-t}}\over
(\t-t)^{n+1\over2}}\2n\int_{\dbR^n}\2n{e^{-\l{|\bar\eta-x|^2\over\t-s}}
\over(\t-s)^{n+1\over2}}d\bar\eta\)d\eta d\t\]\2n\sup_{(\t,\eta)\in[s,T]\times\dbR^n}|\th(\t,\t,\eta,\eta)
\1n-\1n\h\th(\t,\t,\eta,\eta)|\\
\ns\ds\les K\int_s^T{1\over(\t-s)^{1\over2}}\sup_{\xi,x\in\dbR^n}|\th_x
(t,\t,\xi,x)-\h\th_x(t,\t,\xi,x)|d\t\\
\ns\ds\qq+K\[\int_s^T\({1\over(\t-s)^{1\over2}}+{1\over(\t-t)^{1\over2}
(\t-s)^{1\over2}}\)d\t\]\sup_{{(t,\t)\in\h\D[s,T]}\atop{\xi,\eta\in\dbR^n}}
|\th(t,\t,\xi,\eta)-\h\th(t,\t,\xi,\eta)|.\ea$$
Then noting $1<p<2$, by Young's inequality, we have
$$\ba{ll}
\ns\ds\(\int_{t\vee S}^T\sup_{\xi,x\in\dbR^n}|\cS[\,\th\,]_x(t,s,\xi,x)
-\cS[\,\h\th\,]_x(t,s,\xi,x)|^pds\)^{1\over p}\\
\ns\ds\les\1n K(T\1n-\1n S)^{1\over2}\1n\[\(\1n\int_{t\vee S}^T\sup_{\xi,x\in\dbR^n}|\th_x
(t,\t,\xi,x)\1n-\1n\h\th_x(t,\t,\xi,x)|^pd\t\)^{1\over p}\2n+\2n\sup_{{(t,s)\in\h\D[S,T]}\atop{\xi,x\in\dbR^n}}|\th(t,s,\xi,x)
-\h\th(t,s,\xi,x)|\]\\
\ns\ds\les K(T-S)^{1\over2}\|\th-\h\th\|_{\sY[S,T]}.\ea$$
Combining the above, we obtain
\bel{}\|\cS[\,\th\,]-\cS[\,\h\th\,]\|_{\sY[S,T]}\les K(T-S)^{p-1\over p}\|\th-\h\th\|_{\sY[S,T]}.\ee
By taking $\h\th=0$, we see that
$$\|\cS[\,\th\,]\|_{\sY[S,T]}\les\|\cS[0]\|_{\sY[S,T]}
+K(T-S)^{p-1\over p}\|\th\|_{\sY[S,T]},\qq\forall\th\in\sY[S,T].$$
Consequently, $\cS:\sY[S,T]\to\sY[S,T]$ and it is a contraction when $\d=T-S>0$ is small. Hence, it admits a unique fixed point on $[T-\d,T]$, which gives a solution to \rf{Th=int} on $\h\D[T-\d,T]\times\dbR^{2n}$. By repeating the same argument, we will be able to get a unique solution of \rf{Th=int} on $\D[0,T]\times\dbR^{2n}$. Then we define $\G$ by \rf{G=int}. This gives the existence of a mild solution $(\G,\Th)$ of \rf{PDE-[0,T]}.

\ms

The argument used to established the contractiveness of the solution map $\cS$ also gives the uniqueness of the mild solution. \endpf

\section{Concluding Remarks}

In this paper, we have derived the representations of adapted solutions of Type-I BSVIEs and adapted M-solutions of Type-II BSVIEs in terms of the solution to forward SDEs via the solutions of representation PDEs. For Type-I BSVIEs, the well-posedness of representation PDE is established in the classical sense, and for Type-II BSVIEs, the well-posedness of representation PDE is established in the mild solution. It remains open at the moment whether the representation PDE for Type-II BSVIEs admits a unique classical solution, which we believe it to be true, under certain conditions.

\ms

On the other hand, our results could also be regarded as Feynman-Kac formula, from which the solutions to the PDE systems of forms \rf{PDE0} and \rf{PDE**-M} can be represented by the solutions to the corresponding BSVIEs.

\ms

It is worthy of pointing out that, to our best knowledge, representation PDEs of form \rf{PDE0} appeared the first time in the study of time-inconsistent optimal control problems (\cite{Yong 2012}, see also \cite{Wei-Yong-Yu 2017, Mei-Yong 2017}). This indicates that there should be some intrinsic relationship between BSVIEs and time-inconsistent optimal control problems. We hope to explore that in our future publications.

\baselineskip 18pt
\renewcommand{\baselinestretch}{1.2}


\begin{thebibliography}{99}


\bibitem{Antonelli 1993} F.~Antonelli, \it Backward-forward stochastic differential equations, \sl Ann. Appl. Probab., \rm 3 (1993), 777--793.

\bibitem{Bellman-Kalaba-Wing 1960} R.~Bellman, R.~Kalaba and G.~M.~Wing, \it Invariant imbedding and the reduction of two-point boundary value problems to initial value problems, \sl Proc. Nat. Acad. Sci. U.S.A., \rm 46 (1960), 1646--1649.

\bibitem{Bellman-Wing 1975} R.~Bellman and G.~M.~Wing, \sl An Introduction to Invariant Imbedding, \rm Wiley-Interscience, New York, 1975.

\bibitem{Bismut 1973} J.-M.~Bismut, \it Th\'eorie Probabiliste du Contr\^ole des Diffusions, \sl Mem. Amer. Math. Soc. \rm 176, Providence, Rhode Island, 1973.

\bibitem{Briand-Delyon-Hu-Pardoux-Stoica 2003} Ph.~Priand, D.~Delyon, Y.~Hu, E.~Pardoux and L.~Stoica, \it $L^p$ solutions of backward stochastic differential equations, \sl Stoch. Proc. Appl., \rm 108 (2003), 109--129.

\bibitem{Briand-Hu 2006} Ph.~Briand and Y.~Hu, \it BSDE with quadratic growth and unbounded terminal value, \sl Probab. Theory Rel. Fields, \rm 136 (2006), 604--618.

\bibitem{Delbaen-Hu-Bao 2011} F.~Delbaen, Y.~Hu and X.~Bao, \it Backward SDEs with superquadratic growth, \sl Probab. Theory Rel. Fields, \rm 150 (2011), 145--192.

\bibitem{Douglas-Ma-Protter 1996} J.~Douglas,~Jr., J.~ Ma and P.~Protter, \it Numerical methods for forward-backward stochastic differential equations, \sl Ann. Appl. Probab., \rm 6 (1996), 940--968.

\bibitem{Duffie-Epstein 1992} D.~Duffie and L.~Epstein, \it Stochastic differential utility, \sl Econometrica, \rm 60 (1992), 353--394.

\bibitem{Ekren-Keller-Touzi-Zhang 2014} I.~Ekren, C.~Keller, N.~Touzi and J.~Zhang, \it On viscosity solutions of path dependent PDEs, \sl Ann. Probab., \rm 42 (2014), 204--236.

\bibitem{El Karoui-Peng-Quenez 1997} N.~El Karoui, S.~Peng, and M.~C.~Quenez, \it Backward stochastic differential equation in finance, \sl Math. Finance, \rm 7 (1997), 1--71.

\bibitem{Friedman 1964} A.~Friedman, Partial Differential Equations of Parabolic Type, Prentice-Hall,  Englewood Cliffs, N.J., 1964.

\bibitem{Hu-Ma 2004} Y.~Hu and J.~Ma, \it Nonlinear Feynman-Kac formula and discrete-functional-type BSDEs with continuous coefficients, Stoch. Proc. Appl., \rm 112 (2004), 23--51.

\bibitem{Hu-Peng 1995} Y.~Hu and S.~Peng, \it Solution of forward-backward stochastic differential equations, \sl Probab. Theory Rel. Fields, \rm 103 (1995), 273--283.

\bibitem{Karatzas-Shreve 1988} I.~Karatzas and S.~E.~Shreve, \sl Brownian Motion and Stochastic Calculus, \rm Springer-Verlag, Berlin, 1988.

\bibitem{Kobylanski 2000} M.~Kobylanski, \it Backward stochastic differential equations and partial differential equations with quadratic growth, \sl Ann. Probab., \rm 28 (2000), 558--602.

\bibitem{Ladyzenskaja 1968} O.~A.~Ladyzenskaja, V.~A.~Solonnikov, and N.~N.~Ural'ceva, \sl Linear and Quasi-linear Equations of Parabolic Type, \rm AMS, Providence, R.I., 1968.

\bibitem{Lin 2002} J.~Lin, \it Adapted solution of a backward stochastic nonlinear Volterra integral equation, \sl Stoch. Anal. Appl., \rm 20 (2002), 165--183.

\bibitem{Ma-Protter-Yong 1994} J.~Ma, P.~Protter, and J.~Yong, \it Solving forward-backward stochastic differential equations explicitly --- a four step scheme, \sl Prob. Theory Rel. Fields, \rm 98 (1994), 339--359.

\bibitem{Ma-Wu-Zhang-Zhang 2015} J.~Ma, Z.~Wu, D.~Zhang, J.~Zhang, \it On well-posedness of forward-backward SDEs --- a unified approach, \sl Ann. Appl. Probab., \rm 25 (2015), 2168--2214.

\bibitem{Ma-Yong 1995} J.~Ma and J.~Yong, \it Solvability of forward backward SDEs and the nodal set of Hamilton-Jaccobi-Bellman Equations, \sl Chin. Ann. Math., \rm 16B (1995), 279--298.

\bibitem{Ma-Yong 1999} J.~Ma and J.~Yong, \sl Forward-backward stochastic differential equations and their applications, \rm Lecture Notes in Mathematics, 1702. Springer-Verlag, Berlin, 1999.

\bibitem{Ma-Yong-Zhao 2010} J.~Ma, J.~Yong and Y.~Zhao, \it Four step scheme for general Markovian forward-backward SDEs, \sl J. Syst. Sci. Complex., \rm 23 (2010), 546--571.

\bibitem{Ma-Zhang 2002} J.~Ma and J.~Zhang, \it Representation theorems for backward stochastic differential equations, \sl Ann. Appl. Probab., \rm 12 (2002), 1390--1418.

\bibitem{Mei-Yong 2017} H.~Mei and J.~Yong, \it Equilibrium strategies for time-inconsistent stochastic switching systems, \rm arXiv:1712.09505v1[math.OC]27 Dec 2017.

\bibitem{Pardoux-Peng 1990} E.~Pardoux and S.~Peng, \it Adapted solution of a backward stochastic differential equation, \sl Systems \& Control Lett., \rm 14 (1990), 55--61.

\bibitem{Peng 1991} S.~Peng, \it A nonlinear Feynman-Kac formula and applications, \sl Control Theory, Stochastc Analysis and Applications, S.~Chen and J.~Yong, eds, \rm World Scientic, Singapore, 1991, 173--184.

\bibitem{Peng 2004} S.~Peng, \it Nonlinear expectations, nonlinear evaluations and risk measures, \sl Stochastic Methods in Finance, \rm 165--253, Lecture Notes in Math., 1856, Springer, Berlin, 2004.

\bibitem{Peng 2010} S.~Peng, \it Backward stochastic differential equation, nonlinear expectation and their applications, \sl Proceedings of the International Congress of Mathematicians, \rm Volume I, 393--432, Hindustan Book Agency, New Delhi, 2010.

\bibitem{Soner-Touzi-Zhang 2012} H.~M.~Soner, N.~Touzi and J.~Zhang, \it  Wellposedness of second order backward SDEs, \sl Probab. Theory Rel. Fields, \rm 153 (2012), 149--190.

\bibitem{Shi-Wang-Yong 2013} Y.~Shi, T.~Wang and J.~Yong, \it Mean-field backward stochastic Volterra integral equations, \sl Discrete Continuous Dynamic Systems, \rm 18 (2013), 1929--1967.

\bibitem{Shi-Wany-Yong 2015} Y.~Shi, T.~Wang and J.~Yong, \it Optimal control problems of forward-backward stochastic Volterra integral equations, \sl Math. Control Rel. Fields, \rm 5 (2015), 613--649.

\bibitem{Wang-Yong 2015} T.~Wang and J.~Yong, \it Comparison theorems for backward stochastic integral equations, \sl Stoch. Proc. Appl., \rm 125 (2015), 1756--1798.

\bibitem{Wang 2016} Y.~Wang, \it A numerical scheme for BSVIEs, \rm arXiv:1605.04865v1 [math.NA] 16 May 2016.

\bibitem{Wei-Yong-Yu 2017} Q. Wei, J. Yong and Z. Yu, Time-inconsistent recursive stochastic optimal control problems, \sl SIAM J. Control Optim., \rm 55 (2017), 4156--4201.


\bibitem{Yong 1997} J.~Yong, \it Finding adapted solution of forward-backward stochastic differential equations --- method of continuation, \sl Probab. Theory Rel. Fields, \rm 107 (1997) 537--572.

\bibitem{Yong 2006} J.~Yong, \it Backward stochastic Volterra integral equations and some related problems, \sl Stoch. Proc. Appl, \rm 116 (2006), 779--795.

\bibitem{Yong 2007} J.~Yong, \it Continuous-time dynamic risk measures by backward stochastic Volterra integral equations, Appl. Anal, \rm 86 (2007), 1429--1442.

\bibitem{Yong 2008} J.~Yong, \it Well-posedness and regularity of backward stochastic Volterra integral equations, \sl Probab. Theory Rel. Fields, \rm 142 (2008), 21--77.

\bibitem{Yong 2010} J.~Yong, \it Forward-backward stochastic differential equations with mixed initial-terminal conditions, \sl Trans. Amer. Math. Soc., \rm 362 (2010), 1047--1096.

\bibitem{Yong 2012} J.~Yong, \it Time-inconsistent optimal control problems and the equilibrium HJB equation, \sl Math. Control Rel. Fields, \rm 2 (2012), 271--329.

\bibitem{Yong 2016} J.~Yong, \it Representation of adapted solutions to backward stochastic Volterra integral equations, \sl Scientia Sinica Mathematica, \rm 47 (2017), 305--345 (in Chinese).

\bibitem{Zhang 2004} J.~Zhang, \it A numerical scheme for BSDEs, \sl Ann Appl. Probab., \rm 14 (2004), 459--488.


\end{thebibliography}
\end{document}